\documentclass[10pt]{amsart}

\usepackage{amsmath}
\usepackage{amsthm}
\usepackage{amssymb}
\usepackage{xypic}

\newcommand{\calb}{\mathcal{B}}
\newcommand{\calc}{\mathcal{C}}
\newcommand{\cald}{\mathcal{D}}

\newcommand{\calf}{\mathcal{F}}

\newcommand{\caln}{\mathcal{N}}

\newcommand{\calr}{\mathcal{R}}
\newcommand{\cals}{\mathcal{S}}

\newcommand{\calu}{\mathcal{U}}

\newcommand{\IC}{{\mathbb C}} 
 
\newcommand{\IF}{{\mathbb F}}

\newcommand{\IQ}{{\mathbb Q}} 
\newcommand{\IR}{{\mathbb R}}
 
\newcommand{\IZ}{{\mathbb Z}}

\newcommand{\dom}{{\rm dom}}   
\newcommand{\im}{{\rm im}}
\newcommand{\coker}{{\rm coker}}
\newcommand{\res}{{\rm res}}

\newcommand{\id}{{\rm id}}
\newcommand{\tr}{{\rm tr}}

\newcommand{\colim}{{\rm colim }}

\newcommand{\semidirect}{\rtimes}

\newcommand{\Tor}{{\rm Tor}}

\newcommand{\GL}{{\rm GL}}        

\newcommand{\NZD}{{\rm NZD}}

\newcommand{\cell}{\mathit{cell}}
\newcommand{\ab}{\mathit{ab}}

\theoremstyle{plain}
\newtheorem{theorem}{Theorem}[section]
\newtheorem{lemma}[theorem]{Lemma}
\newtheorem{corollary}[theorem]{Corollary}
\newtheorem{proposition}[theorem]{Proposition}
\newtheorem{conjecture}[theorem]{Conjecture}
\newtheorem{addendum}[theorem]{Addendum}
\newtheorem{consequence}[theorem]{Consequence}
\newtheorem{question}[theorem]{Question}
\newtheorem*{corollaryohne}{Corollary}
\newtheorem*{theoremohne}{Theorem}

\theoremstyle{definition}
\newtheorem{definition}[theorem]{Definition}
\newtheorem{example}[theorem]{Example}

\newtheorem{note}[theorem]{Note}

\theoremstyle{remark}
\newtheorem{remark}[theorem]{Remark}
\newtheorem{digression}[theorem]{Digression}

\date{\today}

\begin{document}

\title{$L^2$-Betti numbers, Isomorphism Conjectures and Noncommutative Localization} 
\author{Holger Reich}
\address{Fachbereich Mathematik\\Universit\"at M\"unster\\ Einsteinstr.~62\\
  48149~M\"unster\\ Germany}
\email{reichh@math.uni-muenster}
\thanks{Research supported by the SFB "Geometrische Strukturen in der
  Mathematik" in M\"unster, \indent Germany.} 

\begin{abstract}
In this paper we want to discuss how the question about the rationality of $L^2$-Betti numbers is related
to the Isomorphism Conjecture in algebraic $K$-theory and why in this context noncommutative localization appears
as an important tool.
\end{abstract}

\maketitle



\tableofcontents


\typeout{---------------------...und los -------------------------}


{\em $L^2$-Betti numbers}\/ are invariants of spaces which are defined analogously
to the ordinary Betti-numbers but they take information about the fundamental group into account
and are a priori real valued.
\vspace{1em}

The {\em Isomorphism Conjecture in algebraic $K$-theory}\/ predicts that $K_0( \IC \Gamma )$, the 
Grothendieck group of finitely generated projective
$\IC \Gamma$-modules, should be computable from the $K$-theory of 
the complex group rings of finite subgroups of $\Gamma$.
\vspace{1em}

Given a commutative ring one can always invert the set of all non-zerodivisors. Elements in
the resulting ring have a nice description in terms of fractions. For noncommutative rings
like group rings this may no longer be the case and other concepts for a  {\em noncommutative localization}\/ can be more suitable for
specific problems.
\vspace{1em}

The question whether $L^2$-Betti numbers are always rational numbers was asked by Atiyah in \cite{Atiyah(1976)}.
The question turns out to be a question about modules over the group ring of the fundamental group $\Gamma$.
In \cite{Linnell(1993)} Linnell was able to answer the question affirmatively if $\Gamma$ belongs to 
a  certain class of groups which contains free groups and is stable under extensions by elementary amenable groups 
(one also needs a bound on the orders of finite subgroups). In fact Linnell proves the stronger result
that there exists a semisimple subring in $\calu \Gamma$, the algebra of operators affiliated to the group von Neumann algebra,
which contains the complex group ring.

The main purpose of this short survey  is to give a conceptional framework for Linnell's result, to explain 
how the question about the rationality of $L^2$-Betti numbers relates to the Isomorphism conjecture,
and why this may involve studying noncommutative localizations of group rings. (The impatient reader should right 
away take a look at Proposition~\ref{strategy}, Theorem~\ref{main-linnell} and Addendum~\ref{addendum}.)

Since probably not every reader is familiar with all three circles of ideas -- $L^2$-Betti numbers -- Isomorphism Conjectures
-- Noncommutative Localization -- the paper contains introductions to all of these.

After a brief introduction to group von Neumann algebras and the notion of $\Gamma$-dimension we proceed to explain
the algebra $\calu \Gamma$ of operators affiliated to a group von Neumann algebra and introduce $L^2$-Betti numbers in a very algebraic
fashion. (Once $\calu \Gamma$ has been defined there is no more need for Hilbert-spaces.) 
Section~\ref{section-atiyah} explains the Atiyah Conjecture and contains in particular 
Proposition~\ref{strategy} which is a kind of strategy
for its proof. That Proposition says that if one can factorize the inclusion $\IC \Gamma \subset \calu \Gamma$ 
over a ring $\cals \Gamma$ with good ring-theoretical properties in such way that a 
certain $K$-theoretic condition is satisfied, then the Atiyah conjecture follows. 
In Section~\ref{section-candidates} we present a number of candidates for the ring $\cals \Gamma$.
To do this we first 
review a number of concepts from the theory of noncommutative localization in Section~\ref{section-localization}.
One of the candidates is 
the universal localization of $\IC \Gamma$ with respect to all matrices that become invertible 
over $\calu \Gamma$. Section~\ref{section-linnell} contains Linnell's result. We would like to emphasize that
the intermediate rings Linnell constructs can also be viewed as universal localizations, see Addendum~\ref{addendum}~(U).
In Section~\ref{section-isomorphism} we discuss the Isomorphism Conjecture which seems to be closely related to the 
$K$-theoretical condition mentioned above.
In the last Section we discuss to what extent the functor $- \otimes_{\IC \Gamma} \calu \Gamma$, which plays an important role 
when one studies $L^2$-Betti numbers, is exact.

The only new result in this paper is the following.
\begin{theoremohne}
Let $\cald \Gamma$ denote the division closure of $\IC \Gamma$ inside $\calu \Gamma$.
The flat dimension of $\cald \Gamma$ over $\IC \Gamma$ is smaller than $1$ for 
groups in Linnell's class which have a bound on the orders of finite subgroups.
\end{theoremohne}
For the division closure see Definition~\ref{definition-division-rational}, 
for Linnell's class of groups see Definition~\ref{definition-linnells-class}.
The result is proven as Theorem~\ref{application} below. As an immediate corollary one obtains.
\begin{corollaryohne} \label{corollary-euler-intro}
If the infinite group $\Gamma$ belongs to Linnell's class $\calc$ and has a bound on the orders of finite subgroups 
then the $L^2$-Euler characteristic (and hence the ordinary one whenever defined) satisfies
\[
\chi^{(2)} ( \Gamma ) \leq 0.
\]
\end{corollaryohne}
We also would like to mention that the above theorem leads to interesting non-trivial examples of stably flat
universal localizations which appear in \cite{Neeman-Ranicki(2002)}.

A reader who is interested in more information about $L^2$-Betti numbers and the Atiyah Conjecture should 
consult the book \cite{Lueck(2002)}. Almost all topics discussed here are also treated there in detail.
More information and further results about the Atiyah Conjecture can be found in 
\cite{Linnell(1998)}, \cite{Grigorchuk-Linnell-Schick-Zuk(2000)}, 
\cite{Schick(2000c)} and \cite{Schick(2002a)}.
\vspace{1em}

{\bf Acknowledgement}

I would like to thank Wolfgang L{\"u}ck who as my thesis advisor introduced me to the Atiyah Conjecture
and the Isomorphism Conjecture. I would also like to thank Thomas Schick. In discussions with him the idea that the 
Corollary above  should be true evolved. Furthermore I would like to thank Andrew Ranicki
and the ICMS in Edinburgh for organizing the lively and interesting Workshop on Noncommutative Localization.

\section{The von Neumann Dimension} \label{section-dimension}
In this section we want to introduce group von Neumann algebras and 
explain a notion of dimension for finitely generated projective modules over such algebras.

For a (discrete) group  $\Gamma$ we denote by $\IC \Gamma$ the complex group ring and by $l^2 \Gamma$
the complex Hilbert space with orthonormal basis $\Gamma$. Each group element operates from the left on
$l^2 \Gamma$. Linearly extending this action we obtain an inclusion
\[
\IC \Gamma \to \calb ( l^2 \Gamma )
\]
into the algebra $\calb ( l^2 \Gamma )$ of bounded linear operators on the Hilbert space $l^2 \Gamma$.
The group von Neumann algebra $\caln \Gamma$ is defined as the closure of $\IC \Gamma$ inside $\calb ( l^2 \Gamma )$
with respect to the weak (or strong, it doesn't matter) operator topology. This algebra is closed under taking the adjoint, i.e.\
$a^{\ast} \in \caln \Gamma$ for every operator $a \in \caln \Gamma$.
\begin{digression}
A von Neumann algebra is by definition a $\ast$-closed subalgebra of the algebra of bounded
linear operators on some Hilbert-space which is closed with respect to the strong (or weak, it doesn't matter) operator
topology. Similarly a  $C^{\ast}$-algebra can be defined as a $\ast$-closed subalgebra of the algebra of bounded operators on some 
Hilbert space which is
closed with respect to the topology given by the operator-norm. 
Every von Neumann algebra is in particular a $C^{\ast}$-algebra.
In the situation described above the 
operator-norm closure of $\IC \Gamma$ inside $\calb ( l^2 \Gamma )$ defines the so called reduced $C^{\ast}$-algebra $C^{\ast}_r \Gamma$
and we  have a natural inclusion of $C^{\ast}_r \Gamma$ in  $\caln \Gamma$. 
\end{digression}
The bicommutant theorem of von Neumann (see for example Theorem~5.3.1 in \cite{Kadison-Ringrose(1983)}) is a first hint that 
the definition of $\caln \Gamma$ is very natural also from a purely algebraic point of view
(at least if we agree to consider $\calb ( l^2 \Gamma )$ as something natural). 
It says that the von Neumann algebra is the double commutant of $\IC \Gamma$, i.e.\
\[
\caln \Gamma = \IC \Gamma^{\prime \prime},
\]
where for a subset $A \subset \calb ( l^2 \Gamma )$ we write 
$A^{\prime}= \{ b \in \calb( l^2 \Gamma ) | ba = ab \mbox{ for all } a \in A \}$ for the commutant of $A$ in $\calb ( l^2 \Gamma )$.

The group von Neumann algebra comes equipped with a natural trace. This trace is given as follows:
\begin{eqnarray*} 
\tr_{\Gamma} : \caln \Gamma & \to & \IC \\
a & \mapsto & \langle a (e ) , e \rangle.
\end{eqnarray*}
Here $\langle - , - \rangle$ denotes the inner product in $l^2 \Gamma$ and $e$ is the unit element of the group considered
as a vector in $l^2 \Gamma$. Applied to an element $a=\sum a_g g$ in the group ring $\IC \Gamma$ the trace yields
the coefficient of the identity element $a_e$. 
Of course we have the trace property $\tr_{\Gamma}( ab ) = \tr_{\Gamma} ( ba )$. \label{argument ???}
Once we have such a trace there is a standard procedure to assign a complex number to each finitely generated projective 
$\caln \Gamma$-module: if $p=p^2=(p_{ij}) \in M_n ( \caln \Gamma )$ is an idempotent matrix over $\caln \Gamma$ which
represents $P$, i.e.\ such that
$P \cong \im( p: \caln \Gamma^n \to \caln \Gamma^n)$ then we set
\begin{eqnarray} \label{definition-dim}
\dim_{ \Gamma} P = \Sigma_{i = 1}^n \tr_{\Gamma} ( p_{ii} )
\end{eqnarray}
and call it the $ \Gamma$-dimension of $P$. We have the following standard facts.
\begin{proposition}
The $ \Gamma$-dimension has the following properties.
\begin{enumerate}
\item
$\dim_{ \Gamma} P$ is a nonnegative real number.
\item
$\dim_{ \Gamma} P$ depends only on the isomorphism class of $P$.
\item {\bf Normalization.}
We have $\dim_{ \Gamma} \caln \Gamma = 1$.
\item {\bf Additivity.}
If $0 \to L \to M \to N \to 0$ is a short exact sequence of finitely generated projective modules
then 
\[
\dim_{ \Gamma} M = \dim_{ \Gamma} L + \dim_{ \Gamma} N.
\]
\item {\bf Faithfulness.}
$\dim_{ \Gamma} P = 0 $ if and only if $P=0$.
\end{enumerate}
\end{proposition}
\begin{proof}
(i) follows since one can always arrange that the idempotent $p=p^2$ in (\ref{definition-dim})
is a projection, i.e.\ $p=p^2=p^{\ast}$ (see for example Proposition~4.6.2 on p.23 in \cite{Blackadar(1986)}).
(v) follows from the fact that the trace is faithful, i.e.\ $\tr( a^{\ast} a) =0$ implies $a = 0$.
(ii)-(iv) are straightforward.
\end{proof}

Let $K_0 ( \caln \Gamma )$ denote the Grothendieck-group of finitely generated projective $\caln \Gamma$-modules
then because of (i)-(iv) above we obtain a homomorphism
\[
\xymatrix{
K_0 ( \caln \Gamma ) \ar[rr]^-{\dim_{ \Gamma }} & & \IR.
         }
\]
We recall some terminology, compare page~5 in \cite{Schofield(1985)}.
\begin{definition}
A projective rank function $\rho$ on a ring $R$ is a homomorphism $\rho: K_0 ( R ) \to \IR$
satisfying $\rho( [R^1] )=1$ and $\rho( [P] ) \geq 0$ for every finitely generated projective
$R$-module $P$. It is called faithful if moreover $\rho( [P] ) =0$ implies $P=0$.
\end{definition}
In this  terminology we can summarize the content of the proposition above by saying that 
$\dim_{ \Gamma} : K_0 ( \caln \Gamma ) \to \IR$ is a faithful projective rank function.

Other natural examples of faithful projective rank functions occur as follows:
Suppose the ring $R$ is embedded in a simple artinian ring $M_n ( D )$, where $D$ is a skew field.
Then $P \mapsto \frac{1}{n} \dim_D P \otimes_R M_n ( D )$ defines a faithful projective rank function on $R$.

We would like to emphasize the following additional properties of the 
$\Gamma$-dimension for $\caln \Gamma$-modules which are  not true for arbitrary projective rank functions.
They give further justification for the the use of the word ``dimension'' in this context.

\begin{proposition} The $\Gamma$-dimension satisfies:
\begin{enumerate}
\item[(v)]
{\bf Monotony.}
The $\caln \Gamma$-dimension is monotone, i.e.\ $P \subset Q$ implies
that $\dim_{ \Gamma} P \leq \dim_{ \Gamma} Q$. 
\item[(vi)]
{\bf Cofinality.}
If $P = \bigcup_{i \in I } P_i$ is a directed union of submodules then 
\[
\dim_{ \Gamma}  P  = \sup_{i \in I } \dim_{ \Gamma } P_i .
\]
\end{enumerate}
\end{proposition}

Of course cofinality implies monotony.
To convince the reader that these properties are not automatic for projective rank functions we would 
like to treat an example.

\begin{example}  \label{example-monotone-false}
Let $\Gamma$ be a free group on two  generators $x$ and $y$. 
By work of Cohn \cite{Cohn(1964)} 
we know that $\IC \Gamma$ is a free ideal ring. 
In particular every finitely generated projective module is free and taking its rank yields an isomorphism
\[
\xymatrix{
K_0 ( \IC \Gamma ) \ar[r]^-{\cong} & \IZ.
         }
\]
This is a faithful projective rank function with values in $\IZ$.
However there is an exact sequence
\[
\xymatrix{
0 \ar[r] &  \IC \Gamma^2 \ar[rr]^-{(x-1 , y-1 )} & &  \IC \Gamma \ar[r] &  \IC \ar[r] &  0
         } 
\]
which shows that the rank function is not monotone.
(Geometrically the above resolution of $\IC$ is 
obtained as the cellular chain complex with complex coefficients of the universal cover $E \Gamma$ of 
the model for the classifying space $B \Gamma$ given by the wedge of two circles.)
\end{example}

In fact one  can always compose $\dim_{\Gamma}$ with the natural map
$K_0 ( \IC \Gamma ) \to K_0 ( \caln \Gamma )$.  In this way we obtain naturally
a faithful projective rank function on $\IC \Gamma$ for every group $\Gamma$. One  rediscovers the example above in the 
case where $\Gamma$ is  the free group on two generators.

\section{The Algebra of Operators affiliated to $\caln \Gamma$.} \label{section-affiliated}

The category of finitely generated projective $\caln \Gamma$-modules has one drawback:
it is not abelian. In particular if we start out with a complex of finitely generated
projective $\caln \Gamma$-modules then the homology modules are not necessarily finitely generated
projective and hence the $\caln \Gamma$-dimension as explained above is a priori not available.
But this is exactly what we would like to do in order to define $L^2$-Betti numbers, i.e.\ we want to consider 
\[
C^{\cell}_{\ast} ( \widetilde{ X } ) \otimes_{ \IZ \Gamma } \caln \Gamma
\]
the cellular chain-complex of the universal covering of a CW-complex $X$ of finite type
tensored up to $\caln \Gamma$ and assign a dimension to the homology modules.

There are several ways to get around this problem.
The traditional way to deal with it is to work with certain Hilbert spaces with an isometric $\Gamma$-operation 
instead of modules, e.g.\
with $l^2 \Gamma^n$ instead of $\caln \Gamma^n$. These Hilbert spaces have a $\Gamma$-dimension
and one (re-)defines the homology as the kernel of the  differentials modulo the {\em closure}\/ of their  images.
This is then again a Hilbert space with an isometric $\Gamma$-action and has a well defined $\Gamma$-dimension.

A different approach is taken 
in \cite{Lueck(1997a)}: {\em finitely presented} $\caln \Gamma$-modules do form an abelian
category (because $\caln \Gamma$ is a semihereditary ring) 
and the $\caln \Gamma$-dimension can be extended to these modules in such a way that the properties
(i)-(vi) still hold. (In fact in \cite{Lueck(1998a)} the $\Gamma$-dimension is even extended to arbitrary
$\caln \Gamma$-modules.)

A third possible approach is to introduce the algebra $\calu \Gamma$ of operators affiliated to $\caln \Gamma$.
This algebra has better ring-theoretic properties and indeed  finitely generated projective $\calu \Gamma$-modules 
do form an abelian  category. Moreover the notion of $\Gamma$-dimension extends to that algebra.
We want to explain this approach in some detail in this section.

Recall that an unbounded operator $a: \dom (a) \to H$ on a Hilbert space $H$ is a linear map which is defined
on a linear subspace $\dom (a) \subset H$ called the domain of $a$. It is called densely defined if $\dom (a )$ is a 
dense subspace of $H$ and it is called closed if its graph considered as a subspace of $H \oplus H$ is closed. 
Each bounded operator is closed and densely defined. For unbounded operators $a$ and $b$ 
the symbol $a \subset b$ means that restricted to the possibly smaller domain of $a$ the two operators coincide.
The following definition goes back to \cite{Murray-Neumann(1936)}.
\begin{definition}[Affiliated Operators]
A closed and densely defined (possibly unbounded) operator $a: \dom (a) \to l^2 \Gamma$ is affiliated to
$\caln \Gamma$ if $ba \subset ab$ for all $b \in \caln \Gamma' $.
The set 
\[
\calu \Gamma = \{ a: \dom (a ) \to l^2 \Gamma  \; | \; 
a \mbox{ is \begin{tabular}{l}  closed, \\ densely defined \\ and affiliated to $\caln \Gamma$ 
\end{tabular} } \} 
\]
is called the algebra of operators affiliated to $\caln \Gamma$.
\end{definition}
\begin{remark}
Each group element $\gamma \in \Gamma$ acts by right multiplication on $l^2 \Gamma$. 
This defines an element $r_{\gamma} \in \caln \Gamma'$ 
(we had $\Gamma$ acting from the left when we defined $\caln \Gamma$).
In order to prove that a closed densely defined operator $a$ is affiliated it suffices to check
that its domain $\dom ( a )$ is $\Gamma$-invariant and that for all vectors $v \in \dom (a )$ we
have  $r_{\gamma} a (v)  = a r_{\gamma} (v )$ for all $\gamma \in \Gamma$. In this sense the affiliated operators
are precisely the $\Gamma$-equivariant unbounded operators. 
\end{remark}
Observe that the naive composition of two unbounded operators $c$ and $d$ yields an operator 
$dc$ which is only defined on $c^{-1}( \dom (d) )$. Similarly addition is only defined on the intersection
of the domains. It is hence not obvious that $\calu \Gamma$ is an algebra.

\begin{proposition}
The set $\calu \Gamma$ becomes a $\IC$-algebra if we define addition and a product as the 
closure of the naive addition respectively composition of operators.
\end{proposition}
\begin{proof}
This is proven in Chapter XVI in \cite{Murray-Neumann(1936)}. A  proof is reproduced
in Appendix~I in \cite{Reich(1999)} and also in Chapter~8 of \cite{Lueck(2002)}.
\end{proof}
The subalgebra of all bounded operators in $\calu \Gamma$ is $\caln \Gamma$. 
In contrast to $\caln \Gamma$ there seems to be no useful topology on $\calu \Gamma$. So we left the
realm of $C^{\ast}$-algebras and $C^{\ast}$-algebraic methods. 
The reason $\calu \Gamma$ is nevertheless very useful is that we have gained
good ringtheoretical properties.
Let us recall the definition of von Neumann regularity. 

\begin{definition} \label{vNregular-definition}
A ring $R$ is called von Neumann regular if one of the following equivalent conditions is satisfied.
\begin{enumerate}
\item
Every $R$-module $M$ is flat, i.e.\ for every module $M$ the functor $- \otimes_R M$ is exact.
\item
Every finitely presented $R$-module is already finitely generated projective.
\item
The category of finitely generated projective $R$-modules is abelian.
\item \label{vNregular-definition-vier}
For all $x \in R$ there exists a $y \in R$ such that $xyx=x$.
\end{enumerate}
\end{definition}
\begin{proof}
For (i) $\Leftrightarrow$ (iv) see for example Theorem~4.2.9 in~\cite{Weibel(1994)}. 
(i) $\Rightarrow$ (ii) follows since every finitely presented flat $R$-module is projective, see Theorem~3.2.7
in \cite{Weibel(1994)}. Since the tensor product is compatible with colimits, directed colimits are exact
and every module is a directed colimit of finitely presented modules
we obtain (ii) $\Rightarrow$ (i). For (ii) $\Rightarrow$ (iii)
one needs to check that cokernels and kernels between finitely generated projectives are again finitely generated
projective. But a cokernel is essentially a finitely presented module. The argument for the kernel and (iii) 
$\Rightarrow$ (ii) are elementary.
\end{proof}

Note that in particular fields, skew fields, simple artinian rings and semisimple rings are von Neumann regular
(every module is projective over such rings).
The first condition says that von Neumann regular rings form a very natural class of rings from
a homological algebra point of view: they constitute  precisely the rings of weak homological dimension $0$.
The last condition, which seems less conceptional to modern eyes,
was von Neumann's original definition \cite{vonNeumann(1936)} and has the advantage that one  can explicitely 
verify it in the case we are interested in. More information about von Neumann regular rings can be found in \cite{Goodearl(1979)}.

\begin{proposition}
The algebra $\calu \Gamma$ is a von Neumann regular ring.
\end{proposition}
\begin{proof}
Using the polar decomposition and functional calculus 
one can explicitely construct a $y$ as it is required in the characterization \ref{vNregular-definition}~(iii) 
of von Neumann regularity given above. Compare Proposition~2.1~(v) in \cite{Reich(2001)}.
\end{proof}

In order to define $L^2$-Betti numbers it remains to establish a notion of dimension for finitely generated projective
$\calu \Gamma$-modules. \label{Ueberleitung O.K?}

\begin{proposition} \label{NGUG-properties}
We have the following facts about the inclusion $\caln \Gamma \subset \calu \Gamma$.
\begin{enumerate}
\item
The natural map $K_0 ( \caln \Gamma ) \to K_0 ( \calu \Gamma )$ is an isomorphism.
In particular there is a $\Gamma$-dimension for finitely generated projective
$\calu \Gamma$-modules which we simply define via the following diagram:
\[
\xymatrix{
K_0 ( \caln \Gamma ) \ar[dr]^{\dim_{\Gamma}} \ar[rr]^-{\cong} & & K_0 ( \calu \Gamma ) \ar[dl]_{\dim_{\Gamma}} \\
& \IR &
         }
\]
\item
The ring $\calu \Gamma$ is the Ore-localization (compare~Proposition~\ref{proposition-ore-localization}) 
of $\caln \Gamma$ with respect to the multiplicative subset
of all non-zerodivisors. In particular $- \otimes_{\caln \Gamma} \calu \Gamma$ is an exact functor.
\end{enumerate}
\end{proposition}
\begin{proof}
See \cite{Reich(2001)} Proposition~6.1~(i) and Proposition~2.1~(iii).
\end{proof}

If we now start with a finitely presented (as opposed to finitely generated projective) 
$\IC \Gamma$-module $M$ then because of \ref{vNregular-definition}~(ii) we know that $M \otimes_{\IC \Gamma} \calu \Gamma$
is a finitely generated projective $\calu \Gamma$-module and it makes sense to consider its $\Gamma$-dimension.

\begin{remark}
The assignment $M \mapsto \dim_{ \Gamma} ( M \otimes_{\IC \Gamma} \calu \Gamma )$
is a Sylvester module rank function for finitely presented $\IC \Gamma$-modules in the sense of Chapter~7 in \cite{Schofield(1985)}.
\end{remark}

We are now prepared to give a definition of $L^2$-Betti numbers using the $\Gamma$-dimension for $\calu \Gamma$-modules.
Let $X$ be a CW-complex of finite type, i.e.\ there are only finitely many cells in each dimension. 
Let $\widetilde{ X }$ denote the universal covering. It carries a natural CW-structure
and a cellular free $\Gamma= \pi_1 ( X )$-action. There is one $\Gamma$-orbit of cells in $\widetilde{X}$
for each cell in $X$ and in particular the cellular chain complex
$C_{\ast}^{\cell} ( \widetilde{X} )$ is a complex of finitely generated free $\IZ \Gamma$-modules.
\begin{definition} \label{definition-l2-betti}
For a CW-complex $X$ of finite type we define its $L^2$-Betti numbers as
\[
b_p^{(2)} ( X ) = \dim_{\calu \Gamma} H_p ( C_{\ast}^{\cell} ( \widetilde{X} ) \otimes_{\IZ \Gamma} \calu \Gamma ).
\]
\end{definition}
Note that by \ref{vNregular-definition}~(iii) the homology modules are finitely generated projective $\calu \Gamma$-modules
and hence have a well defined $\calu \Gamma$-dimension.

\begin{remark} \label{remark-dimension-for-arbitrary}
As already mentioned it is possible to extend the notion of $\Gamma$-dimension to  arbitrary
$\caln \Gamma$-modules in such a way that one still has ``additivity'' and 
``cofinality'' \cite{Lueck(1998a)}. Of course one has to allow the value $\infty$, and in cases where this value
occurs one has to interpret ``additivity'' and ``cofinality'' suitably. In \cite{Reich(2001)} it is shown that analogously
there is a $\Gamma$-dimension for arbitrary $\calu \Gamma$-modules which is compatible with the one for 
$\caln \Gamma$-modules in the sense that
for an $\caln \Gamma$-module $M$ we have 
\begin{eqnarray} \label{dim-compatible}
\dim_{\calu \Gamma} M \otimes_{\caln \Gamma} \calu \Gamma = \dim_{\caln \Gamma } M.
\end{eqnarray}
Both notions of extended dimension can be used to define $L^2$-Betti numbers for arbitrary spaces by working with
the singular instead of the cellular chain complex. 
From \ref{NGUG-properties}~(ii) we conclude that for a complex $C_{\ast}$ of $\caln \Gamma$-modules we have 
\[
H_{\ast} ( C_{\ast} \otimes_{\caln \Gamma} \calu \Gamma ) = H_{\ast} ( C_{\ast} ) \otimes_{\caln \Gamma} \calu \Gamma.
\]
If we combine this with (\ref{dim-compatible}) we see that the two possible definitions of $L^2$-Betti numbers coincide.
In the following we will not deal with $L^2$-Betti numbers in this generality. We restrict our
attention to CW-complexes of finite type and hence to finitely generated projective $\calu \Gamma$-modules.
\end{remark}

In order to illustrate the notions defined so far we would like to go through two  easy examples.
\begin{example} \label{example-finite}
Suppose $\Gamma$ is a finite group of order $\# \Gamma$. In this case all the functional analysis is obsolete.
We have $\IC \Gamma = \caln \Gamma = \calu \Gamma$ and $l^2 \Gamma=\IC \Gamma$. 
A finitely generated projective module $P$ is just a finite dimensional complex $\Gamma$-representation.
One can check that
\[
\dim_{\Gamma} P = \frac{1}{\# \Gamma} \dim_{\IC} P.
\]
\end{example}

\begin{example} \label{example-infinite-cyclic}
Suppose $\Gamma = C$ is the infinite cyclic group written multiplicatively with generator $z \in C$.
In this case (using Fourier transformation) the Hilbert-space $l^2 \Gamma$ can be identified with $L^2 ( S^1 )$, the square integrable
functions on the unit circle equipped with the standard normalized measure $\mu=\frac{1}{2 \pi} dz$. 
Under this correspondence the group element $z$ corresponds to the function $z \mapsto z$,
where we think of $S^1$ as embedded in the complex plane. The algebras $\IC \Gamma$, $\caln \Gamma$ and $\calu \Gamma$
can be identified as follows:
\begin{eqnarray*}
\IC \Gamma & \leftrightarrow & \IC [ z^{\pm 1} ] 
\mbox{ \begin{tabular}{l} Laurent-polynomials considered \\ as functions on $S^1$ \end{tabular} } \\
\caln \Gamma & \leftrightarrow & L^{\infty} ( S^1 )
\mbox{ \begin{tabular}{l} essentially bounded \\ functions on $S^1$ \end{tabular} } \\
\calu \Gamma & \leftrightarrow & L ( S^1 )
\mbox{ \begin{tabular}{l} measurable functions \\ on $S^1$ \end{tabular} } 
\end{eqnarray*}
The action on $L^2 (S^1 )$ in each case is simply given by multiplication of functions. The trace on $\caln \Gamma$ becomes
the integral $f \mapsto \int_{S^1} f d \mu $. 
For a measurable subset $A \subset S^1$ let $\chi_A$ denote its characteristic function. Then $p=\chi_A$ is a projection
and $P_A =p L^{\infty} ( S^1)$ is a typical finitely generated projective $L^{\infty} (S^1)$-module. 
We have
\[
\dim_{\Gamma} P_A = \tr_{\Gamma} ( p ) = \int_{S^1} \chi_A d \mu = \mu ( A ).
\]
In particular we see that every nonnegative real number can occur as the $\Gamma$-dimension of a finitely generated
projective $\caln \Gamma$- or $\calu \Gamma$-module. The module $L^{\infty} ( S^1 ) / (z-1) L^{\infty} ( S^1 )$ is an
example of a module which becomes trivial (and hence projective) over $L(S^1)$, because $(z-1)$ becomes invertible.
In fact one can show that the there is an isomorphism
\[
\xymatrix{
K_0( L^{\infty} (S^1) ) \cong K_0( L ( S^1 ) ) \ar[r]^-{\cong} & L^{\infty}(S^1 ; \IZ ),
         }
\]
where $L^{\infty}(S^1 ; \IZ )$ denotes the space of integer valued measurable bounded functions on $S^1$, compare 
Proposition~6.1~(iv) in \cite{Reich(2001)}.
Every such function can be written in a unique way as a finite sum 
$f = \sum_{n=- \infty}^{\infty} n \cdot \chi_{A_n}$ with $A_n=f^{-1} ( \{ n \} ) \subset S^1$ and corresponds to 
$\sum_{n=-\infty}^{\infty} n  [P_{A_n}]$ under the above isomorphism.
\end{example}

Once we have the  notion of $L^2$-Betti numbers it is natural to define 
\[
\chi^{(2)}( X ) = \sum (-1)^i b_i^{(2)} ( X ).
\]
A standard argument shows that for a finite CW-complex this $L^2$-Euler characteristic coincides with the ordinary Euler-characteristic.
But in fact since $L^2$-Betti numbers tend to vanish more often than the ordinary Betti-numbers the  
$L^2$-Euler characteristic is often defined in cases where the ordinary one is not. We also define
$L^2$-Betti numbers and the $L^2$-Euler characteristic of a group as
\[
b_p^{(2)} ( \Gamma )= b_p^{(2)}( B \Gamma ) \mbox{ and } \chi^{(2)} ( \Gamma )= \chi^{(2)}( B \Gamma ).
\]
As an example of an application we would like to mention the following result
which is due to Cheeger and Gromov \cite{Cheeger-Gromov(1986)}. 
\begin{theorem}
Let $\Gamma$ be a group which contains an infinite amenable normal subgroup, then
\[
b_p^{(2)}( \Gamma ) = 0 \mbox{ for all } p \mbox{, and hence } \chi^{(2)} ( \Gamma ) = 0.
\]
\end{theorem}

\section{The Atiyah Conjecture} \label{section-atiyah}

The question arises which real numbers do actually occur as values of $L^2$-Betti numbers. 
This question was asked by Atiyah in \cite{Atiyah(1976)} where he first introduced the notion of $L^2$-Betti numbers.
(The definition of $L^2$-Betti numbers at that time only applied to manifolds and was given in terms of the 
Laplace operator on the universal covering.)
It turns out that the  question about the values can be phrased as a question about the passage from finitely presented 
$\IZ \Gamma$- or $\IQ \Gamma$-modules to 
$\calu \Gamma$-modules.
\begin{proposition}
Let $\Lambda$ be an additive subgroup of $\IR$ which contains $\IZ$. Let $\Gamma$ be a finitely presented group. 
The following two statements are equivalent.
\begin{enumerate}
\item
For all CW-complexes $X$ of finite type with fundamental group $\Gamma$ and all $p \geq 0$ we have 
\[
b_p^{(2)} ( X ) \in \Lambda.
\]
\item
For all finitely presented $\IZ \Gamma$-modules $M$ we have
\[
\dim_{\calu \Gamma} ( M \otimes_{\IZ \Gamma } \calu \Gamma ) \in \Lambda.
\]
\end{enumerate}
\end{proposition}
\begin{proof}
Using the additivity of the dimension and the fact that the finitely generated free 
modules of the complex $C_{\ast}^{\cell} ( \widetilde{X} ) \otimes_{\IZ \Gamma} \calu \Gamma$
have integer dimensions (ii) $\Rightarrow$ (i) is straightforward. For the reverse direction
one needs to construct a CW-complex $X$ with fundamental group $\Gamma$ 
such that the presentation matrix of $M$ appears as the, say $5$-th differential in $C_{\ast}^{\cell}( \widetilde{X} )$
whereas the $4$-th differential is zero. This is possible by standard techniques.
For details see Lemma~10.5 in \cite{Lueck(2002)}.
\end{proof}

More generally one can induce up finitely presented modules over $R \Gamma$ for every coefficient ring $R$ with
$\IZ \subset R \subset \IC$ and ask about the values of the corresponding $\Gamma$-dimensions. 
Let $S \subset R$ be a multiplicatively closed subset. Since each finitely presented $(S^{-1}R) \Gamma$-module
is induced from a finitely presented $R \Gamma$-module (clear denominators in a presentation matrix) 
we can without loss of generality assume that $R$ is a field.
In the following we will work for simplicity with the maximal choice $R = \IC$. 

Let us describe a candidate for $\Lambda$.
We denote by $\frac{1}{\# \calf in \Gamma } \IZ$ the additive subgroup of $\IR$ which is 
generated by the set of numbers $\{ \frac{1}{H} \; | \; H \mbox{ a finite subgroup of } \Gamma \}$.
If there is a bound on the orders of finite subgroups then $\frac{1}{\# \calf in \Gamma} \IZ = \frac{1}{l} \IZ$
where $l$ is the least common multiple of the orders of finite subgroups.
If $\Gamma$ is torsionfree then $\frac{1}{ \# \calf in \Gamma } \IZ = \IZ$.

The following Conjecture turned out to be too optimistic in general
(compare Remark~\ref{counterexample} below). 
But it still has a chance of being true if one additionally assumes a bound on the orders of finite subgroups.
\begin{conjecture}[Strong Atiyah Conjecture] \label{strong-Atiyah}
Let $M$ be a finitely presented $\IC \Gamma$-module then
\[
\dim_{\calu \Gamma} M \otimes_{\IC \Gamma} \calu \Gamma \in \frac{1}{\# \calf in \Gamma} \IZ.
\]
\end{conjecture}

We will see below in \ref{atiyah-implies-zero} that this Conjecture implies the Zero-Divisor Conjecture.

\begin{remark}
As explained above the conjecture makes sense with any field $\IF$ such that $\IQ \subset \IF \subset \IC$
as coefficients for the group ring. With $\IF = \IQ$ the conjecture is equivalent
to the corresponding conjecture about the values of $L^2$-Betti numbers. The conjecture with $\IF= \IC$ clearly implies
the conjecture formulated with smaller fields. 
\end{remark}

To get a first idea let us dicuss the Conjecture in the easy case where $\Gamma$ is the infinite cyclic group.
We have already seen in Example~\ref{example-infinite-cyclic} that in this case the inclusion
$\IC \Gamma \subset \calu \Gamma$ can be identified with $\IC [ z^{\pm 1} ] \subset L( S^1 )$, the Laurent
polynomials considered as functions on $S^1$ inside the algebra of all measurable functions on $S^1$.
Clearly $\IC \Gamma$ corresponds to $\IC [ z^{\pm 1} ]$.
The crucial observation now is that 
in this case we find a field in between $\IC \Gamma$ and $\calu \Gamma$.
Let $\IC( z )$ denote the field of fractions of the polynomial ring $\IC [ z ]$ then we have
\[
\IC [ z^{\pm 1} ] \subset \IC ( z ) \subset L ( S^1 ).
\]
Now let $M$ be a finitely presented $\IC [z^{\pm 1}]$-module then $M \otimes_{\IC [z^{\pm 1} ] } \IC ( z )$ is a finitely
generated free $\IC (z)$-module because $\IC(z)$ is a field 
and hence $M \otimes_{\IC [z^{\pm 1} ] } L (S^1 )$ is a finitely generated free $L(S^1)$-module.
In particular its $\Gamma$-dimension is an integer as predicted by Conjecture~\ref{strong-Atiyah}.

Note that $\IC ( z )$ is not contained in the group von Neumann algebra $L^{\infty} (S^1 )$ because
a rational function like for example $z \mapsto \frac{1}{z-1}$ which has a pole on $S^1$ can not be 
essentially bounded. It hence was crucial for this proof that we had the algebra of affiliated
operators $\calu \Gamma$, here $L( S^1 )$, available.

The following generalizes these simple ideas.

\begin{proposition} \label{strategy}
Suppose the inclusion map $\IC \Gamma \to \calu \Gamma$ factorizes over a ring 
$\cals \Gamma$ such that the following two conditions are fulfilled.
\begin{enumerate}
\item[(K)]
The composite map 
\[
\xymatrix{
\colim_{H \in \calf in \Gamma} K_0 ( \IC H ) \ar[r] & K_0 ( \IC \Gamma )  \ar[r] & K_0 (\cals \Gamma )
         }
\]
is surjective.
\item[(R)]
The ring $\cals \Gamma$ is von Neumann regular.
\end{enumerate}
Then Conjecture~\ref{strong-Atiyah} holds for the group $\Gamma$.
\end{proposition}

In the source of the map in (K) the colimit is taken over the finite subgroups of $\Gamma$.
The structure maps in the colimit are induced by inclusions $K \subset H$ and conjugation maps
$c_g: H \to H^g$, $ h \mapsto ghg^{-1}$.

We will see below (compare Theorem~\ref{main-linnell}) 
that there is a reasonably large class of groups for which a factorization of the inclusion
$\IC \Gamma \to \calu \Gamma$ as required above is known to exist.
In order to prove Proposition~\ref{strategy} we need one more fact about $\Gamma$-dimensions.
\begin{proposition} \label{compatible-with-induction}
The $\Gamma$-dimension is compatible with induction, i.e.\ if $G$ is a subgroup of $\Gamma$ then there is a natural
inclusion $\calu G \subset \calu \Gamma$ and for a finitely generated projective
$\calu G$-module $P$ we have
\[
\dim_{\calu \Gamma} P \otimes_{\calu G} \calu \Gamma = \dim_{\calu G} P.
\]
\end{proposition}
\begin{proof}
There exists a natural inclusion $i: \calu G \to \calu \Gamma$ which extends the 
inclusion $i: \caln G \to \caln \Gamma$ because $\calu G$ is the Ore localization of $\caln G$.
The latter inclusion is compatible with the trace, i.e.\ $\tr_{\Gamma}( i(a) ) = \tr_G ( a )$ for $a \in \caln \Gamma$,
see Lemma~1.24 in \cite{Lueck(2002)}. The claim follows from these facts.
\end{proof}

\begin{proof}[Proof of Proposition~\ref{strategy}]
Let $M$ be a finitely presented $\IC \Gamma$-module. Then also $M \otimes_{\IC \Gamma} \cals \Gamma$ is finitely presented
and hence finitely generated projective by \ref{vNregular-definition}~(ii) because we assume that $\cals \Gamma$ is von Neumann regular.
In particular $M \otimes_{\IC \Gamma} \cals \Gamma$ defines a class in $K_0 ( \cals \Gamma )$. Our second assumption implies
that this class comes from $\colim_{H \in \calf in \Gamma} K_0 ( \IC H )$ via the natural map.
It remains to check that the composition
\[
\xymatrix{
\colim_{H \in \calf in \Gamma} K_0 ( \IC H ) \ar[r] &
K_0 ( \IC \Gamma ) \ar[r] & K_0 ( \cals \Gamma ) \ar[r] &
K_0( \calu \Gamma ) \ar[r]^-{\dim_{\Gamma}} & \IR 
         }
\]
lands inside the subgroup $\frac{1}{\# \calf in \Gamma } \IZ$ of $\IR$. But from Example~\ref{example-finite} together 
with Proposition~\ref{compatible-with-induction} we conclude that for a finite subgroup $H$ and
a finitely generated projective $\IC H$-module $P$ we have
\[
\dim_{\calu \Gamma} P \otimes_{\IC H} \calu \Gamma = \dim_{\IC H} P = \frac{1}{\# H} \dim_{\IC} P.
\]
\end{proof}

\begin{remark} 
From~\ref{vNregular-definition}~\ref{vNregular-definition-vier} it follows that the homomorphic 
image of a von Neumann regular ring is again von Neumann regular. 
In particular the image of $\cals \Gamma$ in $\calu \Gamma$ would be von Neumann regular if $\cals \Gamma$ is.
(But it is not clear that the induced map for $K_0$ is surjective, compare Question~\ref{question-surjective}.)
\end{remark}
\begin{note} \label{note-vNr-bound-semisimple}
Suppose $\cals \Gamma$ is a subring of $\calu \Gamma$ which contains $\IC \Gamma$.
If we assume the properties (K) and (R) and additionally we assume that $\Gamma$ has a bound on the orders of finite subgroups,
then $\cals \Gamma$ is semisimple.
\end{note}
\begin{proof}
The assumptions imply that the projective rank function $P \mapsto \dim_{\Gamma} P \otimes_{\cals \Gamma} \calu \Gamma$
for finitely generated $\cals \Gamma$ modules takes values in $\frac{1}{l} \IZ$, where $l$ is the least common multiple
of the orders of finite subgroups. Since each finitely generated projective $\cals \Gamma$-module is a subset of a $\calu \Gamma$-module
it is easy to see that the projective rank function is faithful. In order to prove that a von Neumann regular ring is semisimple
it suffices to show that there are no infinite chains of ideals, see page 21 in \cite{Goodearl(1979)}.
Since each ideal is a direct summand of
$\cals \Gamma$ and each subideal of a given ideal is a direct summand this can be checked using the faithful projective rank function with values in $\frac{1}{l} \IZ$.
\end{proof}

\begin{remark} \label{counterexample} 
The lamplighter group is the semidirect product of $\IZ$ and $\bigoplus_{-\infty}^{\infty}  \IZ/ 2$ where $\IZ$ acts via shift
on $\bigoplus_{-\infty}^{\infty} \IZ/ 2$. The orders of finite subgroups that occur are precisely all powers of $2$.
Conjecture~\ref{strong-Atiyah} hence predicts $\IZ [ \frac{1}{2} ]$ as the range for the dimensions.
However in \cite{Grigorchuk-Linnell-Schick-Zuk(2000)} a finitely presented $\IQ \Gamma$-module is constructed whose $\Gamma$-dimension
is $\frac{1}{3}$. \label{Further reading...???}
\end{remark}

\section{Noncommutative Localization} \label{section-localization}

Our next aim is to present several candidates for the ring $\cals \Gamma$ which appears in Proposition~\ref{strategy}.
In order to do this we first want to fix some language and review a couple of concepts from the theory of 
localization for noncommutative rings. For more on this subject the reader should consult Chapter II in \cite{Stenstroem(1975)},
Chapter~7 in \cite{Cohn(1985)} and Chapter~4 in \cite{Schofield(1985)}.\\

\subsection*{Ore Localization} 

Classically the starting point for the localization of rings is the wish that certain {\em elements}\/ in the ring
should become invertible. In mathematical terms we have the following universal property.
\begin{definition}
Let $T \subset R$ be a subset which does not contain any zero-divisors. 
A ring homomorphism $f:R \to S$ is called $T$-inverting if $f(t)$ is 
invertible for all $t \in T$. A $T$-inverting ringhomomorphism $i:R \to R_T$ is called universally $T$-inverting
if it has the following universal property: given any $T$-inverting ring homomorphism $f:R \to S$
there exists a unique ring homomorphism $\Phi: R_T \to S$ such that 
\[
\xymatrix{
& & R_T \ar[d]^-{\Phi} \\
R \ar[urr]^-i \ar[rr]^-f & & S
         }
\]
commutes.
\end{definition}

A generator and relation construction shows that there always exists a universal 
$T$-inverting ring and as usual 
it is unique up to canonical isomorphism. Given a ring homomorphism $R \to S$ let us agree to write
\[
T(R \to S)
\]
for the set of elements in $R$ which become invertible in $S$. If one replaces 
$T$ by $\overline{T}=T(R \to R_T)$ the universal inverting ring does not change.
We can hence always assume that $T$ is multiplicatively closed. 
A natural maximal choice for $T$ is 
the set  $\NZD ( R )$ of all non-zerodivisors of $R$.

If the ring $R$ is commutative it is well known that there is a model for $R_T$ whose elements are ``fractions''
or more precisely equivalence classes of pairs $(a,t) \in R \times T$. For noncommutative rings the situation
is more complicated. It goes back to Ore  \label{reference !!!!} that under 
a suitable assumption such a calculus of fractions still exists.
\begin{definition}
A multiplicatively closed subset $T \subset R$ which does not contain zero-divisors or zero
itself satisfies the right Ore-condition if for given $(a,s)\in R \times T$ there always exists a $(b,t) \in R \times T$
such that $at=sb$.
\end{definition}
It is clear that this condition is necessary if a calculus of right fractions exists because we need to be able to 
write a given wrong way (left) fraction $s^{-1}a$ as $bt^{-1}$.
It is a bit surprising that this is the only condition.
\begin{proposition} \label{proposition-ore-localization}
Let $T \subset R$ be a multiplicatively closed subset without zero divisors which satisfies the right Ore
condition, then there exists a ring $RT^{-1}$ and a universal $T$-inverting ring-homomorphism
$i: R \to R T^{-1}$ such that every element of $RT^{-1}$ can be written as $i(a)i(t)^{-1}$ with
$(a,t) \in R \times T$.
\end{proposition}
\begin{proof}
Elements in $RT^{-1}$ are equivalence classes of pairs $(a,t) \in R \times T$. The pair $(a,t)$ is 
equivalent to $(b,s)$ if there exist elements $u$, $v \in R$ such that $au=bv$, $su=tv$ and $su=tv \in S$.
For more details see Chapter II in \cite{Stenstroem(1975)}.
\end{proof}

\begin{remark} \label{remark-ore-flat}
Ore-localization is an exact functor, i.e.\ $RT^{-1}$ is a flat $R$-module, see page~57 in \cite{Stenstroem(1975)}.
\end{remark}

\begin{example} \label{example-ore-false}
Let $\Gamma$ be the free group on two generators $x$ and $y$.
The group ring $\IC \Gamma$ does not satisfy the Ore condition with respect to the set $\NZD ( \IC \Gamma )$ of
all non-zerodivisors. Let $C \subset \Gamma$ be the infinite cyclic subgroup generated by $x$. Now $x-1$ is 
a non-zerodivisor since it becomes invertible in $\calu C$ (compare Example~\ref{example-infinite-cyclic})
and therefore in the overring $\calu \Gamma$. In fact every non-trivial element in $\IC \Gamma$ is a non-zerodivisor
since one can embed $\IC \Gamma$ in a skew field. The Ore condition would imply the existence of $(b,t) \in \IC \Gamma \times
\NZD ( \IC \Gamma )$ with $(y-1)t = (x-1) b$ alias
\[
(x-1)^{-1} (y-1) = bt^{-1} .
\]
This implies that $(-b, t)^{tr}$ is in the kernel of the map $(x-1 , y-1 ): \IC \Gamma^2 \to \IC \Gamma$.
But this map is injective, compare Example~\ref{example-monotone-false}.
\end{example}

\subsection*{Localizing Matrices}

Instead of elements one can try to invert maps.
Let $\Sigma$ be a set of homomorphisms between right $R$-modules. A ring homomorphism $R \to S$
is called $\Sigma$-inverting if for every map $\alpha \in \Sigma$ the induced map $\alpha \otimes_R \id_S$
is an isomorphism.

\begin{definition}
A $\Sigma$-inverting ring homomorphism $i: R \to R_{\Sigma}$ is called universal $\Sigma$-inverting if it has the
following universal property. Given any $\Sigma$-inverting ring homomorphism $f:R \to S$ there exists a unique
ring homomorphism $\Psi: R_{\Sigma} \to S$ such that the following diagram commutes.
\[
\xymatrix{
& & R_{\Sigma} \ar[d]^-{\Psi} \\
R \ar[rr]^-f \ar[urr]^-{i} &  & S .
         }
\]
\end{definition}

From now on let us assume that $\Sigma$ is a set of matrices over $R$. 
For a ring homomorphism $R \to S$ we will write
\[
\Sigma ( R \to S )
\]
for the set of all matrices over $R$ which become invertible over $S$.  One can always replace
a given set of matrices $\Sigma$ by $\overline{\Sigma} = \Sigma ( R \to R_{\Sigma} )$ without changing the universal 
$\Sigma$-inverting ringhomomorphism.
There are different constructions
which prove the existence of a universal $\Sigma$-inverting ring homomorphism. 
One possibility is a generator and relation construction where one starts with the free ring on a set of symbols
$\overline{a}_{i,j}$ where $(a_{i,j})$ runs through the matrices in $\Sigma$ and imposes the 
relations which are given in matrix form as $\overline{A} A = A \overline{ A } = 1$, compare
Theorem~2.1 in \cite{Cohn(1985)}.
For more information the reader should consult 
Chapter~7 in \cite{Cohn(1985)} and Chapter~4 in~\cite{Schofield(1985)}.

Another construction due to Malcolmson~\cite{Malcolmson(1982)} (see also \cite{Beachy(1993)}), \label{in this volume?}
a kind of calculus of fractions for matrices, 
allows a certain amount of control over the ring $R_{\Sigma}$. 

As an easy example we would like to mention the following:
A set of matrices is lower multiplicatively closed if $1 \in \Sigma$ and $a$, $b \in \Sigma$ implies that
\[
\left( \begin{array}{cc} a & 0 \\ c & b \end{array} \right) \in \Sigma
\]
for arbitrary matrices $c$ of suitable size. Observe that $\Sigma( R \to S)$ is always lower multiplicatively closed.

\begin{proposition}[Cramer's rule] \label{cramers-rule}
Let $R$ be a ring and $\Sigma$ be a lower multiplicatively closed 
set of matrices over $R$ then every matrix $a$ over $R_{\Sigma}$ satisfies an equation
of the form
\[
s \left( \begin{array}{cc} 1 & 0 \\ 0 & a \end{array} \right) 
  \left( \begin{array}{cc} 1 & x \\ 0 & 1 \end{array} \right) = b
\]
with $s \in \Sigma$, $x \in  M ( R_{\Sigma} )$ and $b \in M ( R )$. 
\end{proposition}
\begin{proof} See Theorem~4.3 on page~53 in \cite{Schofield(1985)}.
\end{proof}
In particular every matrix $a$ over $R_{\Sigma}$ is stably associated over $R_{\Sigma}$ to a matrix $b$
over $R$, i.e.\ there exist invertible matrices $c$, $d \in GL(R_{\Sigma})$ such that
\[
c \left( \begin{array}{cc} a & 0 \\ 0 & 1_n \end{array} \right) d^{-1} = 
\left( \begin{array}{cc} b & 0 \\ 0 & 1_m  \end{array} \right)
\]
with suitable $m$ and $n$.\\

\subsection*{Division Closure and Rational Closure}

Recall that for a given ring homomorphism $R \to S$ we denoted by $T( R \to S)$ the set of all elements
in $R$ which become invertible in $S$ and by $\Sigma( R \to S)$ the set of all matrices over $R$ that become
invertible over $S$.
The universal localizations $R_{T(R \to S)}$ and $R_{\Sigma( R \to S )}$ come with a natural map to $S$.
In the case where $R \to S$ is injective one may wonder whether these maps embed the universal localizations into
$S$. 
The intermediate rings in the following definition
serve as potential candidates for such embedded versions of the universal localizations.

\begin{definition} \label{definition-division-rational}
Let $S$ be a ring.
\begin{enumerate}
\item
A subring $R \subset S$ is called division closed in $S$ if $T(R \subset S) = R^{\times}$, i.e.\ for 
every element $r \in R$ which is invertible in $S$ the inverse $r^{-1}$ lies already in $R$.
\item
A subring $R \subset S$ is called rationally closed in $S$ if $\Sigma( R \subset S)= \GL (R)$, i.e.\
for every matrix $A$ over $R$ which is invertible over $S$ the entries of the inverse matrix $A^{-1}$
are all in $R$.
\item
Given a subring $R \subset S$ the division closure of $R$ in $S$ denoted 
\[
\cald ( R \subset S)
\]
is the smallest division closed subring of $S$ which contains $R$.
\item
Given a subring $R \subset S$ the rational closure of $R$ in $S$ denoted by 
\[
\calr ( R \subset S)
\]
 is the smallest rationally closed subring of $S$ containing $R$.
\end{enumerate}
\end{definition}
Note that the intersection of division closed intermediate rings is again division closed and similarly for rationally closed
rings. This proves the existence of the division and rational closure.
Moreover we really have closure-operations, i.e.
\begin{eqnarray*}
\cald( \cald(R \subset S) \subset S) & = & \cald( R \subset S) \mbox{ and } \\
\calr( \calr(R \subset S) \subset S) & = & \calr( R \subset S).
\end{eqnarray*}
In \cite[Chapter~7, Theorem~1.2]{Cohn(1985)} it is shown that the set
\begin{eqnarray} \label{cohns-description}
\{ a_{i,j} \in S \; | \; ( a_{i,j} ) \mbox{ invertible over } S, \; ( a_{i,j} )^{-1} \mbox{ matrix over } R \}
\end{eqnarray}
is a subring of $S$ and that it is rationally closed. Since this ring is contained in $\calr( R \subset S)$
the two rings coincide.
The following observation is very useful in our context.
\begin{proposition} \label{universally-closed}
A von Neumann regular ring $R$ is division closed and rationally closed in every overring.
\end{proposition}
\begin{proof}
Suppose $a \in R$ is not invertible in $R$, then the corresponding multiplication map $l_a: R \to R$
is not an isomorphism. Therefore the kernel or the cokernel is non-trivial. Both split of as direct summands
because the ring is von Neumann regular. The corresponding projection onto the kernel or cokernel is given
by left multiplication with a suitable idempotent. This idempotent shows that $a$ must be a zerodivisor and 
hence can not become invertible in any overring. A matrix ring over a von Neumann regular ring is again
von Neumann regular and the same reasoning applied to matrix rings over $R$ yields that $R$ is also rationally
closed in every overring.
\end{proof}

In particular note that once we know that the division closure $\cald( R \subset S)$
is von Neumann regular then it coincides with the rational closure $\calr ( R \subset S )$.
The following proposition relates the division respectively rational closure to the 
universal localizations $R_{T(R \subset S)}$ and $R_{\Sigma ( R \subset S )}$.

\begin{proposition} \label{abstract-embedded}
Let $R \subset S$ be a ring extension.
\begin{enumerate}
\item
The map $R_{T( R \subset S ) } \to S$ given by the universal property factorizes
over the division closure.
\[
\xymatrix{ 
& R_{T ( R \subset S )}  \ar[d]^-{\Phi} &  \\
R \ar[r]^-{\subset}  \ar[ur] & \cald ( R \subset S) \ar[r]^-{\subset} & S
         }
\]
\item
If the pair $(R , T(R \subset S))$ satisfies the right Ore condition, then $\Phi$ is an isomorphism.
\item
The map $R_{\Sigma ( R \subset S )} \to S$ given by the universal property factorizes over the rational closure.
\[
\xymatrix{
& R_{\Sigma ( R \subset S )}  \ar@{->>}[d]^-{\Psi} & \\
R \ar[ur] \ar[r]^-{\subset}  & \calr ( R \subset S ) \ar[r]^-{\subset} & S
         }
\]
The map $\Psi$ is always surjective.
\end{enumerate}
\end{proposition}
\begin{proof}
(i) This follows from the definitions. (ii) Note that $T( R \subset S)$ always consists of 
non-zerodivisors. Thus we can choose a ring of right fractions as a model for $R_{T(R \subset S )}$. 
Every element in $\im \Phi$ is of the form $at^{-1}$ with $t \in T( R \subset S)$. Such an element is invertible 
in $S$ if and only if $a \in T ( R \subset S)$. We see that the image
of $\Phi$ is division closed and hence $\Phi$ is surjective. On the other hand the abstract fraction $at^{-1} \in 
R T(R \subset T)^{-1}$ is zero if and only if $a =0$ because $T( R \subset S)$ contains no zerodivisors, so $\Phi$ is injecitve.
(iii) Only the last statement is maybe not obvious. By Cohn's description of the rational closure (compare~(\ref{cohns-description}))
we need to find a preimage for $a_{i,j}$, where $(a_{i,j})$ is a matrix invertible over $S$ whose inverse lies over $R$.
The generator and relation construction of the universal localization immediately gives such an element.
\end{proof}

In general it is not true that the map $\Psi$ is injective.

\section{Some Candidates for $\cals \Gamma$}  \label{section-candidates}

We are now prepared to describe the candidates for the ring $\cals \Gamma$ which appears in Proposition~\ref{strategy}.
We consider the ring extension $\IC \Gamma \subset \calu \Gamma$ and define
\begin{eqnarray*} 
\cald \Gamma  & = & \cald ( \IC \Gamma \subset \calu \Gamma ), \\ 
\calr \Gamma  & = & \calr ( \IC \Gamma \subset \calu \Gamma ), \\
\IC \Gamma_T & = & \IC \Gamma_{T ( \IC \Gamma \subset \calu \Gamma )} \mbox{ and } \\
\IC \Gamma_{\Sigma } & = & \IC \Gamma_{ \Sigma( \IC \Gamma \subset \calu \Gamma ) } .
\end{eqnarray*}

These rings are organized in the following diagram
\begin{eqnarray} \label{diagram-candidates}
\xymatrix{
& \IC \Gamma_T \ar[r] \ar[d] & \IC \Gamma_{\Sigma} \ar[d] & \\
\IC \Gamma \ar[r]^-{\subset} \ar[ur] & \cald \Gamma \ar[r]^-{\subset}  & \calr \Gamma \ar[r]^-{\subset}  & \calu \Gamma 
         }
\end{eqnarray}

A first hint that the rational or division closure may be a good candidate for $\cals \Gamma$ is the following result
which is implicit in~\cite{Linnell(1993)}. At the same time its proof  illustrates
the usefulness of Cramer's rule~\ref{cramers-rule}.

\begin{proposition} \label{RG-skew}
If $\Gamma$ is a torsionfree group then the Strong Atiyah Conjecture~\ref{strong-Atiyah} implies that $\calr \Gamma$ is 
a skew field. 
\end{proposition}
\begin{proof}
For $x \in \calu \Gamma$ let $l_x: \calu \Gamma \to \calu \Gamma$ denote left multiplication with $x$.
From the additivity and faithfulness of the dimension it follows that $x$ is invertible if and only if 
$\dim \im ( l_x ) = 1$ or equivalently $\dim \ker ( l_x )=0$ or equivalently $\dim \coker (l_x) = 0$.
Now let $X$ be a matrix over $\calr \Gamma$ then by \ref{abstract-embedded}~(iii) we know that we can lift it to a matrix
over $\IC \Gamma_{\Sigma}$. Using Cramer's rule~\ref{cramers-rule} and projecting down again we see that we can find invertible matrices
$A$ and $B \in GL ( \calr  \Gamma )$ such that
\[
C= A \left( \begin{array}{cc}  1_n & 0 \\ 0 & X \end{array} \right) B 
\]
is a matrix over $\IC \Gamma$. In particular if $0 \neq x \in \calr  \Gamma $ then for $X=(x)$ we know that 
there exists an $n$ such that
\[
\dim \im (l_x) +n = \dim ( \im \left( \begin{array}{cc} 1_n & 0 \\ 0 & l_x \end{array} \right) ) = \dim ( \im C ) \in \IZ
\]
because we assume for the matrix $C$ over $\IC \Gamma$ that the dimension of its image is an integer. It follows that 
$\dim \im (l_x) = 1$ and hence that $x$ is invertible in $\calr \Gamma$.
\end{proof}

\begin{note} \label{note-DG-RG}
If one of the rings $\cald \Gamma$ or  $\calr \Gamma$ is a skew field then so is the other and the two coincide.
\end{note}
\begin{proof}
If $\cald \Gamma$ is a skew field then it is also rationally closed, see~\ref{universally-closed}.
If $\calr \Gamma$ is a skew field then $\cald \Gamma$ is a division closed subring of a skew field and hence
itself a skew field.
\end{proof}
\begin{corollary} \label{atiyah-implies-zero}
The Atiyah Conjecture~\ref{strong-Atiyah} implies the Zero-Divisor Conjecture, i.e.\ the conjecture
that the complex group ring of a torsionfree group does not contain any zero-divisors.
\end{corollary}

\begin{remark}
One can show that for a torsionfree amenable group the Atiyah Conjecture~\ref{strong-Atiyah} is equivalent to  the 
Zero-Divisor Conjecture, see Lemma~10.16 in \cite{Lueck(2002)}.
\end{remark}

Another natural question is in how far the rings discussed above depend functorially on the group.
Since an arbitrary group homomorphism $G \to G^{\prime}$ does not induce a map from $\calu G$ to $\calu G^{\prime}$
we can not expect functoriality but at least we have the following.
\begin{note}
An injective group homomorphism induces maps on the rings $\cald \Gamma$, $\calr \Gamma$, $\IC \Gamma_{T}$ and 
$\IC \Gamma_{\Sigma}$. These maps are compatible with the maps in diagram~\ref{diagram-candidates}
above.
\end{note}
\begin{proof} 
We already know that the inclusion $\IC \Gamma \subset \calu \Gamma$ is functorial for injective group homomorphisms.
Let $G$ be a subgroup of $\Gamma$. Since $\calu G$ is von Neumann regular it is division closed and rationally closed
in every overring, compare~\ref{universally-closed}. Therefore $\cald \Gamma \cap \calu G$ is division-closed in 
$\calu \Gamma$ and $\cald G \subset \cald \Gamma \cap \calu G \subset \cald \Gamma$.
Analogously one argues for the rational closure. One immediately checks that $T( \IC G \subset \calu G) \subset
T( \IC \Gamma \subset \calu \Gamma)$ and $\Sigma ( \IC G \subset \calu G ) \subset \Sigma ( \IC \Gamma \subset \calu \Gamma )$.
The universal properties imply the statement for $\IC \Gamma_T$ and $\IC \Gamma_{\Sigma}$. 
\end{proof}

\section{Linnell's Result} \label{section-linnell}

Before we state Linnell's result we would like to introduce the class of groups it applies to.
\begin{definition}[Linnell's class of groups] \label{definition-linnells-class}
Let $\calc$ be the smallest class of groups which has the following properties.
\begin{enumerate}
\item[(LC1)]
Free groups are contained in $\calc$.
\item[(LC2)]
If $1 \to G \to \Gamma \to H \to 1$ is an exact sequence of groups such that $G$ lies in $\calc$ and
$H$ is finite or finitely generated abelian then $\Gamma$ lies in $\calc$.
\item[(LC3)]
The class $\calc$ is closed under directed unions, i.e.\ if a group $\Gamma= \bigcup_{i \in I} \Gamma_i$ is a directed union
of subgroups $\Gamma_i$ which lie in $\calc$ then $\Gamma$ also lies in $\calc$.
\end{enumerate}
\end{definition}

To put this definition into perspective we would like to make a couple of remarks.
\begin{remark}
(i) If one replaces (LC1) above by the requirement that the trivial group belongs to $\calc$ one 
obtains the smaller class of elementary amenable groups. Compare  
\cite{Chou(1980)} and \label{.precise.ref..} \cite{Kropholler-Linnell-Moody(1988)}. 
Elementary amenable groups are in particular amenable (see \cite{Day(1957)}) but it is not easy to
find amenable groups that are not elementary amenable \cite{Grigorchuk(1998)}. 
A group which contains a non-abelian free subgroup is not amenable.

(ii) One can show that if $\Gamma$ lies in $\calc$ and $A$ is an elementary amenable normal subgroup then $\Gamma /A$
also belongs to $\calc$.

(iii) The class $\calc$ is closed under free products.
\end{remark}

In \cite{Linnell(1993)} Linnell proves Conjecture~\ref{strong-Atiyah} for groups in the class $\calc$
which additionally have a bound on the orders of finite subgroups. In fact by carefully investigating the
proof given there one can obtain the following statements.
\begin{theorem} \label{main-linnell}
Suppose the group $\Gamma$ lies in $\calc$ and has a bound on the orders of finite subgroups then 
\begin{enumerate}
\item[ (K) ]
The composition 
\[
\colim_{H \in \calf in \Gamma} K_0 ( \IC H ) \to K_0 ( \IC \Gamma ) \to K_0 ( \cald \Gamma )
\]
is surjective.
\item[ (R) ]
The ring $\cald \Gamma$ is semi-simple and hence $\cald \Gamma= \calr \Gamma$ by Proposition~\ref{universally-closed}.
\end{enumerate}
\end{theorem}

As already mentioned this result is essentially contained in \cite{Linnell(1993)}.
In the above formulation it is proven in \cite{Reich(1999)}. The proof is published in Chapter~10
in~\cite{Lueck(2002)}. Below we will only make a couple of remarks about the proof.

Since we formulated the theorem with the division closure $\cald \Gamma$ the reader may get the impression that
this is the best candidate for an intermediate ring $\cals \Gamma$ as in Proposition~\ref{strategy}. But in fact 
the situation is not so clear. 
We already stated that $\cald \Gamma = \calr \Gamma$ when the theorem applies. Moreover one can show the following.
\begin{addendum} \label{addendum} $\mbox{}$
\begin{enumerate}
\item[(U)]
In the situation of Theorem~\ref{main-linnell} 
the natural map $\IC \Gamma_{ \Sigma } \to \calr \Gamma$ is an 
isomorphism and hence $\cald \Gamma = \calr \Gamma \cong \IC \Gamma_{\Sigma}$.
\item[(O)]
If $\Gamma$ lies in the smaller class of elementary amenable groups and has a bound on the orders
of finite subgroups then $\IC \Gamma$ satisfies the the right Ore condition with respect
to the set $\NZD( \IC \Gamma )$ of all non-zerodivisors and this set coincides with 
$T( \IC \Gamma \subset \calu \Gamma)$. Hence $\IC \Gamma_T$ can be realized as a ring of fractions
and the natural map $\IC \Gamma_T \to \cald \Gamma$ is an isomorphism, compare~\ref{abstract-embedded}~(iii).
\end{enumerate}
\end{addendum}
The statement (O) about the Ore localization appears already in \cite{Linnell(1991)}

We will now make some comments about the proof of Theorem~\ref{main-linnell} and the Addendum~\ref{addendum}.
As one might guess from the description of the class of groups to which the Theorem (and the Addendum) applies the proof proceeds via 
transfinite induction on the class of groups, i.e.\ one proves the following statements.
\begin{enumerate}
\item[(I)]
(K), (R) and (U) hold for free groups.
\item[(II)]
If $1 \to G \to \Gamma \to H \to 1$ is  an extension of groups where $H$ is finite or infinite cyclic and 
(K), (R) and (U) hold for $G$ then they hold for $\Gamma$. Similar with (O) replacing (U).
\item[(III)]
If $\Gamma$ is the directed union of the subgroups $\Gamma_i$ and (K), (R) and (U) hold for all $\Gamma_i$ then
they hold for $\Gamma$ if $\Gamma$ has a bound on the orders of finite subgroups. Similar with (O) replacing (U).
\end{enumerate}

\label{induction-principle}
(I) The Kadison Conjecture says that there are no non-trivial idempotents in the group $C^{\ast}$-algebra
of a torsionfree group. Linnell observed that Connes conceptional proof of this Conjecture for the free group on
two generators given in \cite{Connes(1985a)} (see also \cite{Julg-Valette(1984)})
can be used to verify the stronger Conjecture~\ref{strong-Atiyah} 
in this case. Combined with Proposition~\ref{RG-skew} and Note~\ref{note-DG-RG} one concludes that 
$\cald \Gamma = \calr \Gamma$ is a skew field. This yields (R) and also (K) since $K_0$ 
of a skew field is $\IZ$. Every finitely generated free group is a subgroup of the free group on two generators
and every free group is a directed union of finitely generated free subgroups. This is used to pass to arbitrary free groups. 

Some rather non-trivial facts about group rings of  free groups 
(see \cite{Hughes(1970)} and \cite{Lewin(1974)}) are used to verify that $\calr \Gamma$ also coincides with the 
universal localization $\IC \Gamma_{\Sigma}$ and hence to verify (U) in this case. Recall that we saw in Example~\ref{example-ore-false}
that (O) is false for free groups.

(II) For information about crossed products we refer the reader to \cite{Passman(1989)} and to
Digression~\ref{digression-crossed-product} below.
If $1 \to G \to \Gamma \to H \to 1$ is an extension of groups then every set-theoretical section $\mu$ 
of the quotient map $\Gamma \to H$ (we can always assume $\mu(e)=e$ and $\mu(g^{-1})=\mu(g)^{-1}$)
allows to describe the group ring $\IC \Gamma$ 
as a crossed product $\IC G \ast H$ of the ring $\IC G$ with the group $H$.
Similarly crossed products $\cald G \ast H $, $\calr G \ast H$ and $\IC G_{\Sigma} \ast H$ exist and can serve as
intermediate steps when one  tries  to prove the statements (R), (K) and (U) for $\Gamma$.
For example there are natural inclusions
\[
\IC G \ast H \to \cald G \ast H \to \cald \Gamma.
\]
If $H$ is a finite group and $\cald G$ is semisimple then $\cald G \ast H$ is semisimple and coincides with $\cald \Gamma$.
It is relatively easy to verify that $\cald G \ast H$ is noetherian and semiprime if $H$ is an infinite cyclic group and then
Goldie's theorem (a criterion for the existence of an Ore-localization, see Section~9.4 in \cite{Cohn(1991)})
together with results from \cite{Linnell(1992)} are used to verify that $\cald \Gamma$ is the Ore-localization
of $\cald G \ast H$ with respect to the set of all non-zerodivisors. This is roughly the line of argument in order to verify that 
condition (R) survives extensions by finite or infinite cyclic groups. Once we know that $\cald \Gamma$ is an Ore localization
of $\cald G \ast H$ we can combine this with the assumption (which implies that $\IC G_{\Sigma} \ast H \to \cald G \ast H$
is an isomorphism) in order to verify (U). Similarly iterating Ore localizations one obtains that (O) is stable under
extensions with an infinite cyclic group.
Moody's Induction Theorem (see Theorem~\ref{moodys-induction-theorem} below) plays  a crucial role in the argument for (K).
Moreover one has to assume that the class of groups which appears in the induction hypothesis is already closed under 
extensions by finite subgroups. Hence one is forced to start the induction with virtually free groups 
and in particular one has to prove that (K) holds for such groups. For this purpose 
results of Waldhausen \cite{Waldhausen(1978a)} about generalized free products can be used.
Moreover the map induced by $\IC \Gamma \to \IC \Gamma_{\Sigma}$ on $K_0$ needs to be studied, compare 
Question~\ref{question-surjective}. Here it is important to deal with the universal matrix localization.

(III) If $\Gamma$ is the directed union of the subgroups $\Gamma_i$, $i \in I$ then 
$\cald \Gamma$ is the directed union of the subrings $\cald \Gamma_i$ and similar for $\calr \Gamma$. 
A directed union of von Neumann regular rings is again von Neumann regular (use Definition~\ref{vNregular-definition}~(iv))
so $\cald \Gamma$ is at least von Neumann regular
if all the $\cald \Gamma_i$ are semisimple. The fact that $K$-theory 
is compatible with colimits yields that (K) holds for $\Gamma$ if it holds for all the $\Gamma_i$. 
Now the assumption on the bound of the  orders of finite subgroups implies that $\cald \Gamma$ is 
even semisimple by Note~\ref{note-vNr-bound-semisimple}. That (U) and (O) are stable under directed unions
is straightforward.

\begin{digression} \label{digression-crossed-product}
A crossed product $R \ast G = ( S, \mu )$ of the ring $R$ with the group $G$ consists of a ring $S$ which contains
$R$ as a subring together with an injective map $\mu: G \to S^{\times}$ such that the following holds.
\begin{enumerate}
\item
The ring $S$ is a free $R$-module with basis $\mu(G)$.
\item
For every $g \in G$ the conjugation map $c_{\mu(g)}: S \to S$, $\mu(g) s \mu(g)^{-1}$ can be restricted to $R$.
\item
For all $g$, $g^{\prime} \in G$ the element $\tau( g,g^{\prime} )= \mu(g) \mu(g^{\prime}) \mu( gg^{\prime} )^{-1}$
lies in $R^{\times}$.
\end{enumerate}
\end{digression}

\section{The Isomorphism Conjecture in algebraic $K$-theory}  \label{section-isomorphism}

Condition~(K) in Proposition~\ref{strategy} requires that the composite map
\[
\xymatrix{
\colim_{H \in \calf in \Gamma} K_0 ( \IC H ) \ar[r] & K_0 ( \IC \Gamma ) \ar[r] & K_0 ( \cals \Gamma )
         }
\]
is surjective. About the first map in this composition there is the following conjecture.
\begin{conjecture}[Isomorphism Conjecture - Special Case] \label{special-isomorphism-conjecture}
For every group $\Gamma$ the map
\[
\xymatrix{
\colim_{H \in \calf in \Gamma} K_0 ( \IC H ) \ar[r] & K_0 ( \IC \Gamma ) 
         }
\]
is an isomorphism. In particular for a torsionfree group $\Gamma$ we expect 
\[
K_0 ( \IC \Gamma ) \cong \IZ.
\]
\end{conjecture}

To form the colimit we understand $\calf in \Gamma$ as the category whose objects
are the finite subgroups of $\Gamma$ and whose morphisms are generated 
by inclusion maps  $K \subset H$ and conjugation maps $c_g: H \to H^g$, $h \mapsto ghg^{-1}$
with $g \in \Gamma$. 
Observe that in the torsionfree case the colimit reduces to $K_0( \IC ) \cong \IZ$.

In fact Conjecture~\ref{special-isomorphism-conjecture} would be a consequence of a much more general conjecture 
which predicts the whole algebraic
$K$-theory of a group ring $R \Gamma$ in terms of the $K$-theory of the coefficients and homological data about the group. 
This more general conjecture is known as the Farrell-Jones Isomorphism Conjecture for algebraic K-theory \cite{Farrell-Jones(1993a)}.
A precise formulation would require a certain amount of preparation and we refer the reader to
\cite{Davis-Lueck(1998)}, \cite{BFJR} and in particular to \cite{Lueck-Reich(2003)} for more information.

Conjecture~\ref{special-isomorphism-conjecture} would have the following consequence.
\begin{consequence}
For every finitely generated projective $\IC \Gamma$-module $P$ we
have 
\[
\dim_{\Gamma } ( P \otimes_{\IC \Gamma} \calu \Gamma ) \in \frac{1}{\# \calf in \Gamma } \IZ.
\]
\end{consequence}
\begin{proof} Use Example~\ref{example-finite} and Proposition~\ref{compatible-with-induction}.
\end{proof}
So if all finitely presented $\IC \Gamma$-modules were  also finitely generated projective then the Isomorphism
Conjecture would imply the Strong Atiyah Conjecture. But of course this is seldom the case.

Conjecture~\ref{special-isomorphism-conjecture} is true for infinite cyclic groups and products of such by the Bass-Heller-Swan
formula \cite{Bass-Heller-Swan(1964)}.
Cohn's results in  \cite{Cohn(1964)} imply the Conjecture for free groups.
Work of Waldhausen 
\cite{Waldhausen(1978a)} deals with generalized free products and HNN-extensions.
(The reader should consult \cite{Ranicki(2001)} and \cite{Neeman-Ranicki(2002)} for 
a ``noncommutative localization''-perspective on these results.)
A version of the following result plays also an important role in the proof of~\ref{main-linnell}.
\begin{theorem}[Moody's induction theorem - Special case] \label{moodys-induction-theorem}
Let $\Gamma$ be a polycyclic-by-finite group then the map 
\[
\colim_{H \in \calf in \Gamma} K_0 ( \IC H ) \to K_0 ( \IC \Gamma )
\]
is surjective.
\end{theorem}
\begin{proof}
See \cite{Moody(1987)}, \cite{Moody(1989)}, \cite{Cliff-Weiss(1988)}  and Chapter~8 in \cite{Passman(1989)}.
\end{proof}

What happens if we replace the complex coefficients in Conjecture~\ref{special-isomorphism-conjecture} by integral coefficients?
Thinking about the situation for finite groups it is at first glance very surprising that 
for infinite groups there are a lot of cases where there are results about $K_0 ( \IZ \Gamma )$ whereas nothing
is known about $K_0 ( \IC \Gamma )$. See for example~\cite{Farrell-Jones(1993a)}.
The reason is that the elements of algebraic $K$-groups of the {\em integral}\/ group ring  have a topological interpretation.
They occur as obstruction groups in certain topological problems. Many people put a lot of effort into solving
these topological problems and each time this is successful one obtains a result about the algebraic $K$-groups of $\IZ \Gamma$.

However with integral coefficients one does not expect an isomorphism as in Conjecture~\ref{special-isomorphism-conjecture}.
In the case where $\Gamma$ is torsionfree one would still expect $K_0( \IZ \Gamma ) \cong \IZ$, but  in general
so called Nil-groups and also negative $K$-groups should enter in a ``computation'' of $K_0( \IZ \Gamma )$.
Moreover by a result of Swan (see Theorem~8.1 in \cite{Swan(1960a)}) 
the map $K_0 ( \IZ H ) \to K_0 ( \IQ H )$ is almost the trivial
map for a finite group $H$, i.e.\ the map on reduced $K$-groups $\tilde{K}_0 ( \IZ H ) \to \tilde{K}_0 ( \IQ H )$ is trivial.
Summarizing: In general in the square
\[
\xymatrix{
\colim_{H \in \calf in \Gamma} K_0 ( \IZ H ) \ar[d] \ar[r] & K_0 ( \IZ \Gamma ) \ar[d] \\
\colim_{H \in \calf in \Gamma} K_0 ( \IQ H ) \ar[r] & K_0 ( \IQ \Gamma ).
         }
\]
neither the upper horizontal arrow nor the vertical arrows are surjective.  We see that the comparison to the integral group ring
is not very useful for the question we are interested in.

The main techniques to prove results about the $K$-theory of $\IZ \Gamma$ stems from 
``controlled topology''. See  \cite{Quinn(1987)}, \cite{Quinn(1979)}, \cite{Quinn(1982a)}, 
\cite{Farrell(1996)}, \cite{Farrell(Trieste)} and \cite{Jones(Trieste)}.
The set-up has been adapted to a more algebraic setting \cite{Pedersen-Weibel(1989)} and this ``controlled algebra''
(see \cite{Pedersen-Weibel(1989)},\cite{Carlsson(1989)} and \cite{Pedersen(2000)}) was used successfully
to obtain ``lower bounds'' for the $K$-theory of group rings with arbitrary coefficients 
under certain curvature conditions 
on the group \cite{Carlsson-Pedersen(1995a)}.

A result about Conjecture~\ref{special-isomorphism-conjecture}  which uses this ``controlled algebra'' is the following
result from \cite{BFJR}. Recall that a ring is called (right)-regular if it is right noetherian and 
every finitely generated right $R$-module admits finite dimensional projective resolution.
\begin{theorem}
Let $\Gamma$ be the fundamental group of a closed riemannian manifold with strictly negative sectional curvature.
Let $R$ be a regular ring, e.g.\ $R = \IC$ then
\[
K_0( R ) \cong K_0 ( R \Gamma ).
\]
Moreover $K_{-n} ( R \Gamma ) = 0$ and $K_1 ( R \Gamma ) = \Gamma_{\ab} \otimes_{\IZ} K_0 ( R ) \oplus K_1 ( R )$, where
$\Gamma_{\ab}$ denotes the abelianized group.
\end{theorem}

The assumption about $\Gamma$ implies that $\Gamma$ is torsionfree so the above verifies Conjecture~\ref{special-isomorphism-conjecture}.

The author is optimistic that in the near future techniques similar to the ones used in \cite{BFJR} 
will lead to further results about Conjecture~\ref{special-isomorphism-conjecture}.
In view of condition (K) in Proposition~\ref{strategy} the following seems to be an important question.
\begin{question} \label{question-surjective}
Are the maps
\begin{eqnarray*}
K_0 ( \IC \Gamma ) & \to & K_0 ( \IC \Gamma_{\Sigma} ), \\
K_0 ( \IC \Gamma ) & \to & K_0 ( \calr \Gamma ) \\
\mbox{ or }  \quad 
K_0 ( \IC \Gamma ) & \to & K_0 ( \cald \Gamma ) 
\end{eqnarray*}
surjective?
\end{question}
Note that this is true for groups in Linnell's class $\calc$ with a bound on the orders of finite subgroups by
Theorem~\ref{main-linnell}~(K).

\section{Exactness Properties}

In this section we want to investigate to what extent the functor
$- \otimes_{\IC \Gamma} \calu \Gamma$ and related functors are exact. Recall that this functor is crucial 
for the definition of $L^2$-Betti numbers, compare Definition~\ref{definition-l2-betti}.
\begin{note} \label{note-elementary-amenable-flat}
If $\Gamma$ is elementary amenable and there is a bound on the orders of finite subgroups 
then $- \otimes_{\IC \Gamma} \calu \Gamma$ is exact.
\end{note}
\begin{proof}
From Addendum~\ref{addendum}~(O) we know that for these groups 
$\cald \Gamma$ is an Ore-localization of $\IC \Gamma$. In particular in this case
$- \otimes_{\IC \Gamma} \cald \Gamma$ is exact. Since by Theorem~\ref{main-linnell}~(R)
$\cald \Gamma$ is also semisimple (and hence von Neumann regular) we know that every
module is flat over $\cald \Gamma$. 
\end{proof}

The following tells us that we cannot always have exactness. 
\begin{note} 
Suppose for the infinite group $\Gamma$ the functor $- \otimes_{\IC \Gamma} \calu \Gamma$ is exact, then all 
$L^2$-Betti numbers and also the Euler-characteristic $\chi^{(2)}( \Gamma )$ of the group $\Gamma$ vanishes.
\end{note}
\begin{proof}
Flatness implies 
\[
H_p ( C_{\ast} ( E \Gamma ) \otimes_{\IZ \Gamma} \calu \Gamma ) = Tor_p^{\IZ \Gamma} ( \IZ , \calu \Gamma )
 = Tor_p^{\IZ \Gamma} ( \IZ , \IZ \Gamma ) \otimes_{\IZ \Gamma} \calu \Gamma =0
\]
for $p >0$. Moreover $b_0^{(2)}(\Gamma)=0$ for every infinite group (see Theorem~6.54~(8)~(b) in \cite{Lueck(2002)}). 
\end{proof}

In particular we see that for the free group on two generators we cannot have exactness.
We saw this phenomenon already in Example~\ref{example-monotone-false} because exactness of 
$- \otimes_{\IC \Gamma} \caln \Gamma$ would contradict the monotony of the dimension. (Recall from
Proposition~\ref{NGUG-properties}~(ii) that $-\otimes_{\caln \Gamma} \calu \Gamma$ is always exact.)

More generally we have.
\begin{note} \label{note-free-group-not-flat}
If $\Gamma$ contains a nonabelian free group, then neither $\cald \Gamma$ nor 
$\calr \Gamma$, $\IC \Gamma_{\Sigma}$ or $\calu \Gamma$ can be flat over $\IC \Gamma$.
\end{note}
\begin{proof}
Every free group contains a free group on two generators. Let $G \subset \Gamma$ be a free subgroup on two generators.
Let $\IC G^2 \to \IC G $ be the injective homomorphism from Example~\ref{example-monotone-false}.
Since $\IC \Gamma$ is flat over $\IC G$ we obtain an injective map $\IC \Gamma^2 \to \IC \Gamma$. On the other hand 
since $\cald G$ is a skew-field we know that the non-trivial kernel of the corresponding map $\cald G^2 \to \cald G$
(which must appear for dimension reasons since $-\otimes_{\cald G} \calu G$ is exact and the $\Gamma$-dimension
is faithful) is a one-dimensional free module
which splits off $\cald G^2$ as a direct summand. The same remains true for every overring of $\cald G$. In particular
for $\cald \Gamma$, $\calr \Gamma$  and $\calu \Gamma$. But also for $\IC \Gamma_{\Sigma}$ because 
$\cald G = \calr G \cong \IC G_{\Sigma}$ by Addendum~\ref{addendum}~(O) and since 
there is a natural map $\IC G_{\Sigma} \to \IC \Gamma_{\Sigma}$. 
\end{proof}

In this context we would also like to mention the following result from \cite{Lueck(1998b)}. 
\begin{theorem}
If $\Gamma$ is amenable then $\caln \Gamma$ (and hence $\calu \Gamma$) 
is dimension-flat over $\IC \Gamma$, i.e.\ for $p >0$ and every $\IC \Gamma$-module $M$ 
we have
\[
\dim_{\caln \Gamma} Tor_p^{\IC \Gamma} ( M , \caln \Gamma ) = \dim_{\calu \Gamma} Tor_p^{\IC \Gamma} ( M , \calu \Gamma) = 0.
\]
\end{theorem}
\begin{proof} See \cite{Lueck(1998b)} or Theorem~6.37 on page~259 in \cite{Lueck(2002)} and recall that $\calu \Gamma$ 
is flat over $\caln \Gamma$ and the $\calu \Gamma$-dimension and the $\caln \Gamma$-dimension are compatible.
\end{proof}

Given these facts it is tempting to conjecture that $- \otimes_{\IC \Gamma} \calu \Gamma$ is exact if and only if
$\Gamma$ is amenable. However in \cite{Linnell-Lueck-Schick(2002)} it is shown that the condition about the bound on the 
orders of finite subgroups in Note~\ref{note-elementary-amenable-flat} is necessary. 

\begin{example}
Let $H$ be a nontrivial finite group and let $H \wr \IZ$ denote the 
semidirect product $\bigoplus_{- \infty}^{\infty} H  \semidirect \IZ$, where $\IZ$ is acting via shift 
on the $\bigoplus_{-\infty}^{\infty} H$. Then neither $\cald \Gamma$ nor $\calu \Gamma$ is flat over $\IC \Gamma$
(see Theorem~1 in \cite{Linnell-Lueck-Schick(2002)}).
\end{example}

The main purpose of this section is to prove the following result which measures 
the deviation from exactness for groups in Linnell's class.
\begin{theorem} \label{application}
Let $\Gamma$ be in the class $\calc$ with
a bound on the orders of finite subgroups,
 then
\[
\Tor_p^{\IC \Gamma } ( - ; \cald \Gamma) =0 \hspace{3em} 
\mbox{for all } p \geq  2.
\]
\end{theorem}
Note that for these groups $\cald \Gamma=\calr \Gamma \cong \IC \Gamma_{\Sigma}$ is semisimple and therefore the functor $ - \otimes_{\cald \Gamma} \calu \Gamma$ is exact. The functor $-\otimes_{\IZ \Gamma } \IC \Gamma$ is always exact. Therefore we
obtain the corresponding statements for $\Tor_p^{\IC \Gamma}( - ; \calu \Gamma)$,
$\Tor_p^{\IZ \Gamma} ( - ; \cald \Gamma)$ and $\Tor_p^{\IZ \Gamma} (- ; \calu \Gamma)$.

As an immediate consequence we obtain interesting examples of stably flat universal localizations.
\begin{corollary}
If $\Gamma$ lies in Linnell's class $\calc$ and has a bound on the orders of finite subgroups then 
$\cald \Gamma \cong \IC \Gamma_{\Sigma}$ is stably flat over $\IC \Gamma$, i.e.\ we have 
\[
Tor_p^{\IC \Gamma} ( \cald \Gamma , \cald \Gamma ) = 0 \hspace{3em} \mbox{ for all } p \geq 1.
\]
\end{corollary}
\begin{proof}
We know that $\cald \Gamma \cong \IC \Gamma_{\Sigma}$ is a universal localization of $\IC \Gamma$ and
hence $\IC \Gamma \to \IC \Gamma_{\Sigma}$ is an epimorphism in the category of rings, see page~56 in \cite{Schofield(1985)}.
By Theorem~4.8~b) in \cite{Schofield(1985)} we know that $Tor_1^{\IC \Gamma} ( \cald \Gamma , \cald \Gamma ) = 0$.
For $p \geq 2$ the result follows from Theorem~\ref{application}.
\end{proof}
Recent work of Neeman and Ranicki~\cite{Neeman-Ranicki(2002)} shows that for universal localizations which are 
stably flat there exists a long exact localization sequence which extends Schofield's localization sequence for universal
localizations
(see Theorem~4.12 in \cite{Schofield(1985)}) to the left. 
In the case of Ore-localizations the corresponding sequence was known for a long time, see \cite{Gersten(1974)}, \cite{Grayson(1980)},
\cite{Weibel-Yao(1992)} and \cite{Thomason-Trobaugh(1990)}.
Observe that because of  Note~\ref{note-free-group-not-flat} we know that whenever $\Gamma$ contains a free group $\IC \Gamma_{\Sigma}$
cannot be an Ore-localization.

Here is another consequence of Theorem~\ref{application}.
\begin{corollary} \label{euler char}
If the infinite group $\Gamma$ belongs to $\calc$ and has a bound 
on the orders of finite subgroups, then
\[
\chi^{(2)}(\Gamma) \leq 0.
\] 
\end{corollary}
\begin{proof}
Since the group is infinite we have $b_0^{(2)}( \Gamma )= 0$. Because of 
\[
H_p(  \Gamma ; \calu \Gamma )=\Tor_p^{\IC \Gamma} ( \IC ; \calu \Gamma)=0 \hspace{3em} \mbox{ for all } p \geq 0
\]
we know that $b_1^{(2)}(\Gamma)$ is the only
$L^2$-Betti number which could possibly be nonzero.
\end{proof}
The $L^{(2)}$-Euler characteristic coincides with the
usual Euler-characteristic and the rational Euler-Characteristic of \cite{Wall(1961)} whenever these are defined.

Before we proceed to the proof of Theorem~\ref{application} we would also like to mention  the following
consequences for $L^2$-homology.

\begin{corollary}[Universal Coefficient Theorem] \label{UCT}\index{$L^2$-homology!universal coefficient theorem for}
Let $\Gamma$ be in $\calc$ with a bound on the orders of finite
subgroups. Then 
there is a universal coefficient theorem for $L^2$-homology: 
Let $X$ be a $\Gamma$-space whose isotropy groups are all finite, then there is an
exact sequence
\[
0 \to H_n (X ; \IZ) \otimes_{ \IZ \Gamma } \calu \Gamma \to
H_n^{\Gamma} ( X ; \calu \Gamma ) \to \Tor_1( H_{n-1} (X ; \IZ) ; \calu \Gamma )
\to 0.
\]
\end{corollary}

\begin{proof}
We freely use the dimension theory for arbitrary $\calu \Gamma$-modules, compare Remark~\ref{remark-dimension-for-arbitrary}.
If $X$ has finite isotropy, then the set of singular simplices also
has only finite isotropy groups. If $H$ is a finite subgroup of $\Gamma$,
then $\IC \left[ \Gamma / H \right] \cong \IC \Gamma \otimes_{\IC H} \IC$ is induced from the 
projective $\IC H$-module $\IC$ and therefore projective. We see that
the singular chain complex with complex coefficients
$C_{\ast}= C_{\ast}^{sing} ( X ; \IC )$ is a complex of projective
$\IC \Gamma$-modules. The $E^2$-term of the K\"unneth spectral sequence (compare Theorem 5.6.4 on 
page 143 in \cite{Weibel(1994)}) \index{spectral sequence!K{\"u}nneth}\index{K{\"u}nneth spectral sequence}
\[
E^2_{pq} = \Tor_p^{\IC \Gamma} ( H_q( C_{\ast} ) ; \cald \Gamma ) \Rightarrow
H_{p+q}(C_{\ast} \otimes \cald \Gamma ) = H_{p+q} ( X ; \cald \Gamma)
\]
is concentrated in two columns. The spectral sequence collapses, and we get
exact sequences
\[
0 \to H_n(X; \IC ) \otimes_{\IC \Gamma} \cald \Gamma \to H_n(X ; \cald \Gamma ) \to
\Tor_1^{\IC \Gamma } ( H_{n-1} (X ; \IC ) ; \cald \Gamma )  \to 0.
\]
Applying the exact functor $- \otimes_{\cald \Gamma } \calu \Gamma$ yields the 
result.
\end{proof}

The proof of Theorem~\ref{application} depends on the following Lemma.

\begin{lemma} \label{Tor lemma}
\begin{enumerate}
\item \label{Tl1}
Let $R \ast G \subset S \ast G $ be compatible with the crossed product structure.
Let $M$ be an $R \ast G$-module. There is a natural isomorphism
of right $S$-modules
\[
\Tor_p^{R \ast G} ( M ; S \ast G) \cong \Tor_p^{R} ( \res^{R \ast G}_R M ; S )
\]
for all $p \geq 0$.
\item \label{Tl2}
Suppose $R \subset S$ is a ring extension and  $R = \bigcup_{i \in I} R_i$ 
is the directed union of the subrings $R_i$. Let $M$ be an 
$R$-module. Then
there is a natural isomorphism of right $S$-modules
\[
\Tor_p^{R}( M ; S ) \cong \colim_{i \in I }\Tor_p^{R_i}( \res^{R}_{R_i} M ; S_i ) \otimes_{S_i} S
\]
for all $p\geq 0$.
\end{enumerate}
\end{lemma}
\begin{proof}
\ref{Tl1} We start with the case $p=0$. We denote the crossed product structure
map by $\mu$, compare~Digression~\ref{digression-crossed-product}. Define a map
\[
h_M: \res^{R \ast G}_R M \otimes_R S \to M \otimes_{R \ast G} S \ast G
\]
by $m \otimes s \mapsto m \otimes s$. Obviously $h$ is a natural transformation
from the functor $\res^{R \ast G}_R (-) \otimes_R S$ to $- \otimes_{R \ast G} S \ast G$.
If $M = R \ast G $ the map $h_{R \ast G}^{-1}: R \ast G \otimes_{R \ast G} S \ast G
\cong S \ast G \to \res^{R \ast G}_R R \ast G \otimes_R S $ given by
$s\mu(g) \mapsto g \otimes c_g^{-1}(s)$ is a well-defined inverse. Since
$h$ is compatible with direct sums we see that $h_{F}$ is an isomorphism
for all free modules $F$.
Now if $M$ is an arbitrary module choose a free resolution $F_{\ast} \to M$ of 
$M$ and apply both functors to 
\[
F_1 \to F_0 \to M \to 0 \to 0.
\]
Both functors are right exact, therefore an application of the five lemma 
yields the result for $p=0$. Now let $P_{\ast} \to M$ be a projective resolution of $M$, then
\begin{eqnarray*}
\Tor_p^{R}( \res_R^{R \ast G} M ; S ) &  = & H_p( \res_R^{R \ast G} P_{\ast} \otimes_R S ) \\
& \stackrel{\cong}{\to} & H_p( P_{\ast} \otimes_{R \ast G} S \ast G ) \\
& = & \Tor_p^R ( M ; S \ast G ) .
\end{eqnarray*}
\ref{Tl2} Again we start with the case $p=0$. 
The natural surjections $\res^R_{R_i} M \otimes_{R_i} S \to M \otimes_R S$ induce
a surjective map 
\[
h_M:\colim_{i \in I} \res^R_{R_i} M \otimes_{R_i} S \to M \otimes_R S
\]
which is natural in $M$. Suppose the element of the colimit represented by $\sum_k m_k \otimes s_k \in 
\res^R_{R_i} M \otimes_{R_i} S$ is mapped to zero in $M \otimes_R S$. By construction
the tensor product $M \otimes_R S$ is the quotient of the free module on the set
$M \times S$ by a relation submodule. But every relation involves only finitely
many elements of $R$, so we can find a $j \in I$ such that $\sum_k m_k \otimes s_k=0$
already in $\res^R_{R_j} M \otimes_{R_j} S$. We see that $h_M$ is an isomorphism.
Now let $P_{\ast} \to M$ be a projective resolution. Since the colimit is an exact functor
it commutes with homology and we get
\begin{eqnarray*}
\colim_{i \in I} \Tor_p^{R_i} ( M ; S) &  = &
\colim_{i \in I} H_p ( \res^R_{R_i} P_{\ast} \otimes_{R_i} S ) \\
& = & H_p( \colim_{i \in I} ( \res_{R_i}^R P_{\ast} \otimes_{R_i} S )) \\
& \stackrel{\cong}{\to } & H_p( P_{\ast} \otimes_{R} S ) \\
& = & \Tor_p^R( M ; S ) .
\end{eqnarray*}
\end{proof}

\begin{proof}[Proof of Theorem~\ref{application}]
The proof works via transfinite induction over the group as for the proof of Linnell's Theorem~\ref{main-linnell} itself, compare
(I), (II) and (III) on page~\pageref{induction-principle}.

{\rm \bf (I)} The statement for free groups is well known: let $\Gamma $ be the free group generated
by the set $S$. The cellular chain complex of the universal covering of the obvious 
$1$-dimensional classifying space gives a projective resolution of the trivial module 
of length one
\[
0 \to \bigoplus_{S} \IC \Gamma \to \IC \Gamma \to \IC \to 0 .
\]
Now if $M$ is an arbitrary $\IC \Gamma$-module we apply $-\otimes_{\IC} M$ to the above 
complex and get a projective resolution of length $1$ for $M$ (diagonal action). (Use that for $P$ a projective $\IC \Gamma$-module 
$P \otimes_{\IC} M$ with the diagonal respectively the left $\Gamma$-action are noncanonically isomorphic $\IC \Gamma$-modules.)

{\rm \bf (II)} The next step is to prove that the statement remains true under extensions by finite 
groups. So let $1 \to G \to  \Gamma \to H \to 1$ be an exact sequence with $H$ finite.
We know that $\cald \Gamma = \cald G \ast H$, see Lemma~10.59 on page 399 in \cite{Lueck(2002)} or Proposition~8.13 in \cite{Reich(1999)}. 
Let $M$
be a $\IC \Gamma$-module, then with Lemma~\ref{Tor lemma} and the induction hypothesis
we conclude 
\begin{eqnarray*}
\Tor_p^{\IC \Gamma}(M; \cald \Gamma) & = & \Tor_p^{\IC G \ast H} ( M ; \cald G \ast H) \\
 & \cong  & 
\Tor_p^{\IC G} ( \res^{\IC G \ast H}_{\IC G} M; \cald G ) \\
& = & 0 \quad \mbox{ for } p>1.
\end{eqnarray*}
The case $H$ infinite cyclic is only slightly more complicated. This time we know from 
Lemma~10.69 in \cite{Lueck(2002)} or Proposition~8.18 in \cite{Reich(1999)} that $\cald \Gamma=(\cald G \ast H) T^{-1}$
is an Ore localization, where $T= T( \cald G \ast H \subset \calu \Gamma)$, i.e.\ the set of all elements in
$\cald G \ast H$ which become invertible in $\calu \Gamma$.
Since Ore localization is an exact functor we get
\begin{eqnarray*}
\Tor_p^{\IC \Gamma}(M ; \cald \Gamma) & = & \Tor_p^{\IC G \ast H} (M ; (\cald G \ast H)T^{-1} ) \\
& \cong & \Tor_p^{\IC G \ast H} (M; \cald G \ast H ) \otimes_{\cald G \ast H } \cald \Gamma  
\end{eqnarray*}
and conclude again with Lemma~\ref{Tor lemma} that this module vanishes if $p>1$.

{\rm \bf (III)} 
The behaviour under directed unions remains to be checked. Let $\Gamma = \bigcup_{i \in I}
\Gamma_i $ be a directed union, then 
using Definition~\ref{vNregular-definition} we see that $\bigcup_{i \in I} \cald \Gamma_i$ is 
von Neumann regular and it is easy to check that it coincides with the division closure $\cald \Gamma$.
Now Lemma~\ref{Tor lemma}
gives
\begin{eqnarray*}
\Tor_p^{\IC \Gamma} ( M ; \cald \Gamma ) & \cong  &
\colim_{i \in I} \Tor_p^{\IC \Gamma_i}( \res^{\IC \Gamma}_{\IC \Gamma_i} M ; \cald \Gamma ) \\
& = & \colim_{i \in I} \Tor_p^{\IC \Gamma_i}( \res^{\IC \Gamma}_{\IC \Gamma_i} M ; \cald \Gamma_i ) \otimes_{\cald \Gamma_i} \cald \Gamma \\
& = & 0 \quad \mbox{ for } p>1.
\end{eqnarray*}
\end{proof}


\typeout{-------------------- References -------------------------------}

\def\cprime{$'$} \def\polhk#1{\setbox0=\hbox{#1}{\ooalign{\hidewidth
  \lower1.5ex\hbox{`}\hidewidth\crcr\unhbox0}}}


\begin{thebibliography}{10}

\bibitem{Atiyah(1976)}
M.~F. Atiyah.
\newblock Elliptic operators, discrete groups and von {N}eumann algebras.
\newblock {\em Ast\'erisque}, 32-33:43--72, 1976.

\bibitem{BFJR}
A.~Bartels, F.~T. Farrell, L.~E. Jones, and H.~Reich.
\newblock On the isomorphism conjecture in algebraic ${K}$-theory.
\newblock Preprintreihe SFB 478 --- Geometrische Strukturen in der Mathematik,
  Heft 175, M\"unster, 2001.

\bibitem{Bass-Heller-Swan(1964)}
H.~Bass, A.~Heller, and R.~G. Swan.
\newblock The {W}hitehead group of a polynomial extension.
\newblock {\em Inst. Hautes \'Etudes Sci. Publ. Math.}, 22:61--79, 1964.

\bibitem{Beachy(1993)}
J.~A. Beachy.
\newblock On universal localization at semiprime {G}oldie ideals.
\newblock In {\em Ring theory (Granville, OH, 1992)}, pages 41--57. World Sci.
  Publishing, River Edge, NJ, 1993.

\bibitem{Blackadar(1986)}
B.~Blackadar.
\newblock {\em ${K}$-theory for operator algebras}.
\newblock Springer-Verlag, New York, 1986.

\bibitem{Carlsson(1989)}
G.~Carlsson.
\newblock Homotopy fixed points in the algebraic {$K$}-theory of certain
  infinite discrete groups.
\newblock In {\em Advances in homotopy theory (Cortona, 1988)}, volume 139 of
  {\em London Math. Soc. Lecture Note Ser.}, pages 5--10. Cambridge Univ.
  Press, Cambridge, 1989.

\bibitem{Carlsson-Pedersen(1995a)}
G.~Carlsson and E.~K. Pedersen.
\newblock Controlled algebra and the {N}ovikov conjectures for ${K}$- and
  ${L}$-theory.
\newblock {\em Topology}, 34(3):731--758, 1995.

\bibitem{Cheeger-Gromov(1986)}
J.~Cheeger and M.~Gromov.
\newblock ${L}\sb 2$-cohomology and group cohomology.
\newblock {\em Topology}, 25(2):189--215, 1986.

\bibitem{Chou(1980)}
C.~Chou.
\newblock Elementary amenable groups.
\newblock {\em Illinois J. Math.}, 24(3):396--407, 1980.

\bibitem{Cliff-Weiss(1988)}
G.~Cliff and A.~Weiss.
\newblock Moody's induction theorem.
\newblock {\em Illinois J. Math.}, 32(3):489--500, 1988.

\bibitem{Cohn(1964)}
P.~M. Cohn.
\newblock Free ideal rings.
\newblock {\em J. Algebra}, 1:47--69, 1964.

\bibitem{Cohn(1985)}
P.~M. Cohn.
\newblock {\em Free rings and their relations}.
\newblock Academic Press Inc. [Harcourt Brace Jovanovich Publishers], London,
  second edition, 1985.

\bibitem{Cohn(1991)}
P.~M. Cohn.
\newblock {\em Algebra. {V}ol. 3}.
\newblock John Wiley \& Sons Ltd., Chichester, second edition, 1991.

\bibitem{Connes(1985a)}
A.~Connes.
\newblock Noncommutative differential geometry.
\newblock {\em Inst. Hautes \'Etudes Sci. Publ. Math.}, (62):257--360, 1985.

\bibitem{Davis-Lueck(1998)}
J.~F. Davis and W.~L{\"u}ck.
\newblock Spaces over a category and assembly maps in isomorphism conjectures
  in ${K}$- and ${L}$-theory.
\newblock {\em $K$-Theory}, 15(3):201--252, 1998.

\bibitem{Day(1957)}
M.~M. Day.
\newblock Amenable semigroups.
\newblock {\em Illinois J. Math.}, 1:509--544, 1957.

\bibitem{Farrell(1996)}
F.~T. Farrell.
\newblock {\em Lectures on surgical methods in rigidity}.
\newblock Published for the Tata Institute of Fundamental Research, Bombay,
  1996.

\bibitem{Farrell-Jones(1993a)}
F.~T. Farrell and L.~E. Jones.
\newblock Isomorphism conjectures in algebraic ${K}$-theory.
\newblock {\em J. Amer. Math. Soc.}, 6(2):249--297, 1993.

\bibitem{Farrell(Trieste)}
T.~Farrell.
\newblock The {B}orel {C}onjecture.
\newblock In {\em Topology of high-dimensional manifolds}, volume 9(1) of {\em
  ICTP Lecture Notes}, pages 225--298. ICTP, Trieste, 2002.

\bibitem{Gersten(1974)}
S.~M. Gersten.
\newblock The localization theorem for projective modules.
\newblock {\em Comm. Algebra}, 2:317--350, 1974.

\bibitem{Goodearl(1979)}
K.~R. Goodearl.
\newblock {\em von {N}eumann regular rings}.
\newblock Pitman (Advanced Publishing Program), Boston, Mass., 1979.

\bibitem{Grayson(1980)}
D.~R. Grayson.
\newblock {$K$}-theory and localization of noncommutative rings.
\newblock {\em J. Pure Appl. Algebra}, 18(2):125--127, 1980.

\bibitem{Grigorchuk(1998)}
R.~I. Grigorchuk.
\newblock An example of a finitely presented amenable group that does not
  belong to the class {E}{G}.
\newblock {\em Mat. Sb.}, 189(1):79--100, 1998.

\bibitem{Grigorchuk-Linnell-Schick-Zuk(2000)}
R.~I. Grigorchuk, P.~A. Linnell, T.~Schick, and A.~{\.Z}uk.
\newblock On a question of {A}tiyah.
\newblock {\em C. R. Acad. Sci. Paris S\'er. I Math.}, 331(9):663--668, 2000.

\bibitem{Hughes(1970)}
I.~Hughes.
\newblock Division rings of fractions for group rings.
\newblock {\em Comm. Pure Appl. Math.}, 23:181--188, 1970.

\bibitem{Jones(Trieste)}
L.~Jones.
\newblock Foliated control theory and its applications.
\newblock In {\em Topology of high-dimensional manifolds}, volume 9(2) of {\em
  ICTP Lecture Notes}, pages 405--460. ICTP, Trieste, 2002.

\bibitem{Julg-Valette(1984)}
P.~Julg and A.~Valette.
\newblock {$K$}-theoretic amenability for {${\rm SL}\sb{2}({\bf Q}\sb{p})$},
  and the action on the associated tree.
\newblock {\em J. Funct. Anal.}, 58(2):194--215, 1984.

\bibitem{Kadison-Ringrose(1983)}
R.~V. Kadison and J.~R. Ringrose.
\newblock {\em Fundamentals of the theory of operator algebras. {V}ol. {I}}.
\newblock Academic Press Inc. [Harcourt Brace Jovanovich Publishers], New York,
  1983.
\newblock Elementary theory.

\bibitem{Kropholler-Linnell-Moody(1988)}
P.~H. Kropholler, P.~A. Linnell, and J.~A. Moody.
\newblock Applications of a new ${K}$-theoretic theorem to soluble group rings.
\newblock {\em Proc. Amer. Math. Soc.}, 104(3):675--684, 1988.

\bibitem{Lewin(1974)}
J.~Lewin.
\newblock Fields of fractions for group algebras of free groups.
\newblock {\em Trans. Amer. Math. Soc.}, 192:339--346, 1974.

\bibitem{Linnell(1991)}
P.~A. Linnell.
\newblock Zero divisors and group von {N}eumann algebras.
\newblock {\em Pacific J. Math.}, 149(2):349--363, 1991.

\bibitem{Linnell(1992)}
P.~A. Linnell.
\newblock Zero divisors and ${L}\sp 2({G})$.
\newblock {\em C. R. Acad. Sci. Paris S\'er. I Math.}, 315(1):49--53, 1992.

\bibitem{Linnell(1993)}
P.~A. Linnell.
\newblock Division rings and group von {N}eumann algebras.
\newblock {\em Forum Math.}, 5(6):561--576, 1993.

\bibitem{Linnell(1998)}
P.~A. Linnell.
\newblock Analytic versions of the zero divisor conjecture.
\newblock In {\em Geometry and cohomology in group theory (Durham, 1994)},
  pages 209--248. Cambridge Univ. Press, Cambridge, 1998.

\bibitem{Linnell-Lueck-Schick(2002)}
P.~A. Linnell, W.~L\"uck, and T.~Schick.
\newblock The {O}re condition, affiliated operators, and the lamplighter group.
\newblock Preprintreihe SFB 478 --- Geometrische Strukturen in der Mathematik,
  Heft 205 M\"unster. To appear in the Proceedings of the school/conference on
  ``High-dimensional Manifold Topology'' in Trieste, May/June 2001, 2002.

\bibitem{Lueck(1997a)}
W.~L{\"u}ck.
\newblock Hilbert modules and modules over finite von {N}eumann algebras and
  applications to ${L}\sp 2$-invariants.
\newblock {\em Math. Ann.}, 309(2):247--285, 1997.

\bibitem{Lueck(1998a)}
W.~L{\"u}ck.
\newblock Dimension theory of arbitrary modules over finite von {N}eumann
  algebras and ${L}\sp 2$-{B}etti numbers. {I}. {F}oundations.
\newblock {\em J. Reine Angew. Math.}, 495:135--162, 1998.

\bibitem{Lueck(1998b)}
W.~L{\"u}ck.
\newblock Dimension theory of arbitrary modules over finite von {N}eumann
  algebras and ${L}\sp 2$-{B}etti numbers. {I}{I}. {A}pplications to
  {G}rothendieck groups, ${L}\sp 2$-{E}uler characteristics and {B}urnside
  groups.
\newblock {\em J. Reine Angew. Math.}, 496:213--236, 1998.

\bibitem{Lueck(2002)}
W.~L\"uck.
\newblock {\em {$L^2$}-invariants: theory and applications to geometry and
  {$K$}-theory}.
\newblock Springer-Verlag, Berlin, 2002.
\newblock Ergebnisse der Mathematik und ihrer Grenzgebiete, 3. Folge, Band 44.

\bibitem{Lueck-Reich(2003)}
W.~L{\"u}ck and H.~Reich.
\newblock The {B}aum-{C}onnes and the {F}arrell-{J}ones conjectures in ${K}$-
  and ${L}$-theory.
\newblock in preparation, 2003.

\bibitem{Malcolmson(1982)}
P.~Malcolmson.
\newblock Construction of universal matrix localization.
\newblock {\em Lecture Notes in Mathematics}, 951, 1982.

\bibitem{Moody(1987)}
J.~A. Moody.
\newblock Induction theorems for infinite groups.
\newblock {\em Bull. Amer. Math. Soc. (N.S.)}, 17(1):113--116, 1987.

\bibitem{Moody(1989)}
J.~A. Moody.
\newblock Brauer induction for ${G}\sb 0$ of certain infinite groups.
\newblock {\em J. Algebra}, 122(1):1--14, 1989.

\bibitem{Murray-Neumann(1936)}
F.~Murray and J.~{von Neumann}.
\newblock On rings of operators.
\newblock {\em Annals of Math.}, 37:116--229, 1936.

\bibitem{Neeman-Ranicki(2002)}
A.~Neeman and A.~Ranicki.
\newblock Noncommutative localization and chain complexes i.
\newblock e-print AT.0109118, 2002.

\bibitem{Passman(1989)}
D.~S. Passman.
\newblock {\em Infinite crossed products}.
\newblock Academic Press Inc., Boston, MA, 1989.

\bibitem{Pedersen(2000)}
E.~K. Pedersen.
\newblock Controlled algebraic {$K$}-theory, a survey.
\newblock In {\em Geometry and topology: Aarhus (1998)}, volume 258 of {\em
  Contemp. Math.}, pages 351--368. Amer. Math. Soc., Providence, RI, 2000.

\bibitem{Pedersen-Weibel(1989)}
E.~K. Pedersen and C.~A. Weibel.
\newblock ${K}$-theory homology of spaces.
\newblock In {\em Algebraic topology (Arcata, CA, 1986)}, pages 346--361.
  Springer-Verlag, Berlin, 1989.

\bibitem{Quinn(1979)}
F.~Quinn.
\newblock Ends of maps. {I}.
\newblock {\em Ann. of Math. (2)}, 110(2):275--331, 1979.

\bibitem{Quinn(1982a)}
F.~Quinn.
\newblock Ends of maps. {I}{I}.
\newblock {\em Invent. Math.}, 68(3):353--424, 1982.

\bibitem{Quinn(1987)}
F.~Quinn.
\newblock Applications of topology with control.
\newblock In {\em Proceedings of the International Congress of Mathematicians,
  Vol. 1, 2 (Berkeley, Calif., 1986)}, pages 598--606, Providence, RI, 1987.
  Amer. Math. Soc.

\bibitem{Ranicki(2001)}
A.~Ranicki.
\newblock Noncommutative localization in topology.
\newblock unpublished notes, 2001.

\bibitem{Reich(1999)}
H.~Reich.
\newblock Group von {N}eumann algebras and related algebras.
\newblock Dissertation Universit\"at G\"ottingen,
  http://www.math.uni-muenster.de/u/lueck/publ/reich/reich.dvi, 1999.

\bibitem{Reich(2001)}
H.~Reich.
\newblock On the {$K$}- and {$L$}-theory of the algebra of operators affiliated
  to a finite von {N}eumann algebra.
\newblock {\em $K$-Theory}, 24(4):303--326, 2001.

\bibitem{Schick(2000c)}
T.~Schick.
\newblock Integrality of ${L}\sp 2$-{B}etti numbers.
\newblock {\em Math. Ann.}, 317(4):727--750, 2000.

\bibitem{Schick(2002a)}
T.~Schick.
\newblock Erratum: ``{I}ntegrality of ${L}\sp 2$-{B}etti numbers'' [{M}ath.
  {A}nn. {\bf 317} (2000), no. 4, 727--750; 1777117].
\newblock {\em Math. Ann.}, 322(2):421--422, 2002.

\bibitem{Schofield(1985)}
A.~H. Schofield.
\newblock {\em Representation of rings over skew fields}.
\newblock Cambridge University Press, Cambridge, 1985.

\bibitem{Stenstroem(1975)}
B.~Stenstr{\"o}m.
\newblock {\em Rings of quotients}.
\newblock Springer-Verlag, New York, 1975.
\newblock Die Grundlehren der Mathematischen Wissenschaften, Band 217, An
  introduction to methods of ring theory.

\bibitem{Swan(1960a)}
R.~G. Swan.
\newblock Induced representations and projective modules.
\newblock {\em Ann. of Math. (2)}, 71:552--578, 1960.

\bibitem{Thomason-Trobaugh(1990)}
R.~W. Thomason and T.~Trobaugh.
\newblock Higher algebraic {$K$}-theory of schemes and of derived categories.
\newblock In {\em The Grothendieck Festschrift, Vol.\ III}, volume~88 of {\em
  Progr. Math.}, pages 247--435. Birkh\"auser Boston, Boston, MA, 1990.

\bibitem{vonNeumann(1936)}
J.~von Neumann.
\newblock {\em Continous geometry}, volume~27 of {\em Princeton Mathematical
  Series}.
\newblock Princeton University Press, Princeton, New Jersey, 1960.

\bibitem{Waldhausen(1978a)}
F.~Waldhausen.
\newblock Algebraic ${K}$-theory of generalized free products. {I}, {I}{I}.
\newblock {\em Ann. of Math. (2)}, 108(1):135--256, 1978.

\bibitem{Wall(1961)}
C.~T.~C. Wall.
\newblock Rational {E}uler characteristics.
\newblock {\em Proc. Cambridge Philos. Soc.}, 57:182--184, 1961.

\bibitem{Weibel-Yao(1992)}
C.~Weibel and D.~Yao.
\newblock Localization for the {$K$}-theory of noncommutative rings.
\newblock In {\em Algebraic $K$-theory, commutative algebra, and algebraic
  geometry (Santa Margherita Ligure, 1989)}, volume 126 of {\em Contemp.
  Math.}, pages 219--230. Amer. Math. Soc., Providence, RI, 1992.

\bibitem{Weibel(1994)}
C.~A. Weibel.
\newblock {\em An introduction to homological algebra}.
\newblock Cambridge University Press, Cambridge, 1994.

\end{thebibliography}

\typeout{------------------------ THE END --------------------------}

\end{document}